\documentclass[11pt]{article}

\usepackage{amssymb,amsthm,amsmath,hyperref}

\numberwithin{equation}{section}

\newtheorem{thm}{Theorem}[section]
\newtheorem{lem}{Lemma}[section]
\newtheorem{cor}{Corollary}[section]
\newtheorem{prop}{Proposition}[section]
\theoremstyle{definition}
\newtheorem{defn}{Definition}[section]
\theoremstyle{remark}
\newtheorem{rem}{Remark}[section]

\allowdisplaybreaks

\begin{document}
\title{Global Behavior of  Spherically Symmetric
Navier-Stokes Equations with Density-Dependent Viscosity\thanks
{This work is supported by NSFC 10571158} }
\author{Ting Zhang\thanks{E-mail:  zhangting79@hotmail.com},  Daoyuan Fang\thanks{E-mail:
dyf@zju.edu.cn}\\
Department of Mathematics, Zhejiang University,\\ Hangzhou 310027,
PR China }
\date{}
\maketitle
\begin{abstract}
In this paper, we study a free boundary problem for compressible
spherically symmetric Navier-Stokes equations without a solid
core. Under certain assumptions imposed on the initial data, we
obtain the global existence and uniqueness of the weak solution,
give some uniform bounds (with respect to time) of the solution
and show that it converges to a stationary one as time tends to
infinity. Moreover, we obtain  the stabilization rate estimates of
exponential type in $L^\infty$-norm and weighted $H^1$-norm of the
solution by constructing some Lyapunov functionals. The results
show that such system is stable under small perturbations, and
could  be applied to the astrophysics.\\
 \textbf{Keywords:} Compressible Navier-Stokes equations; density-dependent
viscosity; free boundary; existence; uniqueness; asymptotic
behavior
\end{abstract}
\section{Introduction.}
We consider the compressible Navier-Stokes equations with
density-dependent viscosity in $\mathbb{R}^n(n\geq2)$, which can
be written in Eulerian coordinates as
    \begin{equation}
      \left \{
\begin{array}{l}
\partial_{\tau}\rho+\nabla\cdot(\rho \vec{u})=0, \\
\partial_{\tau}(\rho \vec{u})+\nabla\cdot(\rho \vec{u}\otimes \vec{u})+\nabla
P=\textrm{div}(\mu(\nabla \vec{u}+\nabla
\vec{u}^{\top}))+\nabla(\lambda \textrm{div}\vec{u})-\rho \vec{f},
\end{array}
      \right.\label{sym-E1.1}
    \end{equation}
Here $\rho$, $P$, $\vec{u}=(\mathrm{u}_1,\ldots,\mathrm{u}_n)$ and
$\vec{f}$ are the density, pressure,  velocity and the external
force, respectively; $\mu=\mu(\rho)$ and $\lambda=\lambda(\rho)$
are two viscosity coefficients.

In this paper,   the initial conditions are
    \begin{equation}
    \rho(\vec{\xi},0)=\rho_0(r),\ r\in[0,b],
    \end{equation}
    \begin{equation}
    \vec{u}(\vec{\xi},0)=u_0(r)\frac{\vec{\xi}}{r},
       r\in(0, b],\ \vec{u}(\vec{\xi},0)|_{\vec{\xi}=0}=u_0(0)=0,
    \end{equation}
where $r=|\vec{\xi}|=\sqrt{\xi^2_1+\cdots+\xi^2_n}$ and $b>0$ is a
constant,
 the  boundary condition is
                \begin{equation}
                 \left\{ (P-\lambda \mathrm{div}\vec{u})\mathrm{Id}-
                 \mu (\nabla \vec{u}+\nabla \vec{u}^{\top})
                  \right\}\cdot \vec{n}=P_{\Gamma}\vec{n},
                  \ \vec{\xi}\in \partial \Omega_\tau\label{sym-E1.4}
                \end{equation}
where $\partial\Omega_\tau=\psi(\partial\Omega_0,\tau)$  is a free
boundary, $\vec{n}$ is the unit outward normal vector of
$\partial\Omega_\tau$ and $P_{\Gamma}>0$ is a external pressure.
Here, $\partial\Omega_0=\{\vec{\xi}\in
\mathbb{R}^n:|\vec{\xi}|=b\}$ is the initial boundary and $\psi$
is the flow of $\vec{u}$:
    \begin{equation}
      \left\{
      \begin{array}{ll}
      \partial_\tau \psi(\vec{\xi},\tau)=\vec{u}(\psi(\vec{\xi},\tau),\tau), &\vec{\xi}\in
      \mathbb{R}^n,\\
      \psi(\vec{\xi},0)=\vec{\xi}.
      \end{array}
      \right.
    \end{equation}

 To simplify the presentation, we only consider the famous
polytropic model, i.e. $P(\rho)=A\rho^\gamma$ with $\gamma>1$ and
$A>0$ being constants.  And we assume that the viscosity
coefficients $\mu$ and $\lambda$ are proportional to
$\rho^\theta$, i.e. $\mu(\rho)=c_1\rho^\theta$ and
$\lambda(\rho)=c_2\rho^\theta$ where $c_1,c_2$ and $\theta$ are
three constants.

For the initial-boundary value
 problem (\ref{sym-E1.1})-(\ref{sym-E1.4}),
     we are  looking for a
spherically symmetric solution $(\rho,\vec{u})$:
    $$
      \rho(\vec{\xi},\tau)=\rho(r,\tau),\ \
    \vec{u}(\vec{\xi},\tau)=u(r,\tau)\frac{\vec{\xi}}{r},
    $$
  with the spherically
symmetric external force
    $$\vec{f}=f(m,r, \tau)\frac{\vec{\xi}}{r},
    \ m(\rho,r)=\int_0^r\rho(s,\tau)s^{n-1}ds,
    \ r>0
    $$
and
 $
 \partial
 \Omega_\tau=\{\vec{\xi}\in\mathbb{R}^n:|\vec{\xi}|=b(\tau),b(0)=b,b'(\tau)=u(b(\tau),\tau)\}.
 $

Then $(\rho,u)(r,\tau)$ is determined by
   \begin{equation}
      \left\{\begin{array}{lll}
      \partial_\tau \rho+\partial_r(\rho u)+\frac{n-1}{r}\rho u&=&0,\\
           \rho(\partial_\tau u+u\partial_ru)
      +\partial_r P&=&(\lambda+2\mu)(\partial^2_{rr}u+\frac{n-1}{r}\partial_ru-\frac{n-1}{r^2}u)
      \\ &&+2\partial_r\mu\partial_ru+\partial_r\lambda(\partial_r u+\frac{n-1}{r}u)-\rho f,
     \end{array}
      \right.\label{sym-E2.1}
    \end{equation}
where $(r,\tau)\in(0,b(\tau))\times(0,\infty)$, with the initial
data
    \begin{equation}
      (\rho,u)|_{\tau=0}=(\rho_0,u_0)(r),
    \  0\leq r\leq b,
    \end{equation}
 the fixed boundary condition
            \begin{equation}
            u|_{r=0}=0,
            \end{equation}
and  the free boundary condition
                \begin{equation}
                 \left.\left\{ P-2\mu \partial_ru-\lambda\left(\partial_ru+\frac{n-1}{r}u\right)
                 \right\}\right|_{r=b(\tau)}
                  =P_\Gamma,
                  \label{sym-E2.4}
                \end{equation}
where  $b(0)=b$, $b'(\tau)=u(b(\tau),\tau)$.

Additionally,  we assume the external force $ f(m, r, \tau)$ and
external pressure $P_\Gamma(\tau)\in C^{1}(\mathbb{R}_+)$
 satisfy
    \begin{equation}
    P_\Gamma(\tau)=P_\infty+\Delta P(\tau),
    \     f(m,r,\tau)=f_{\infty}(m,r)+\Delta f(m,r,\tau),\label{sym-f1}
    \end{equation}
for all $r\geq 0$ and     $\tau\geq0$, with
    \begin{equation}
       f_\infty(m,r)=\frac{Gm}{r^{n-1}},\ m(\rho,r)=\int^r_0\rho s^{n-1}ds,
       \ \Delta f(m,r,\tau)\in
    C^{1}(\mathbb{R}_+\times\mathbb{R}_+\times\mathbb{R}_+)\label{f2}
        \end{equation}
         \begin{equation}
       \|\Delta f(\cdot,\cdot,\tau)\|_{L^\infty(\mathbb{R}_+\times\mathbb{R}_+)}\leq f_1(\tau),
    \          \|(\partial_r\Delta f,\partial_\tau\Delta
          f)(\cdot,\cdot,\tau)\|_{L^\infty(\mathbb{R}_+\times\mathbb{R}_+)}\leq f_2(\tau),\label{f4}
        \end{equation}
        \begin{equation}
      f_1,\Delta P\in L^\infty\cap
      L^1(\mathbb{R}_+),\ (\Delta P)',\ f_2\in {L^2}(\mathbb{R}_+),
        \end{equation}
where $\mathbb{R}_+=[0,\infty)$, $P_\infty$ and  $G$ are two
positive constants, perturbations $(\Delta P,\Delta f)$ tend to
$0$ as $\tau\rightarrow \infty$ in some weak sense. $f_\infty$ is
the precise expression for its own gravitational force and $\Delta
f$ expresses the influence of the outside gravitational force, in
the astrophysical case (with spherical symmetry). $P_\Gamma$ also
could express the influence of the surface tension force on the
free boundary. This system can be treated as a simple model of one
fluid in $\Omega_\tau$, whose evolution is influenced by the
gravitational force and the external pressure generated by the
other substance in $\mathbb{R}^n\backslash \Omega_\tau$. We study
the stabilization problem of such system, which could  be applied
to the astrophysics.

Now, we consider the stationary problem,  namely
    \begin{equation}
    (P(\rho_\infty))_r=-\rho_\infty f_\infty(m(\rho_\infty,r),r)
    \label{E-sym1}
    \end{equation}
in an interval $r\in(0,l_\infty)$ with the end $l_\infty$
satisfying
    \begin{equation}
     P( \rho_\infty(l_\infty))=P_\infty,
    \end{equation}
        \begin{equation}
          \int_0^{l_\infty}\rho_\infty r^{n-1}
          dr=M:=\int_0^{b}\rho_0 r^{n-1}dr.
          \label{E-sym3}
        \end{equation}
The unknown quantities are the stationary density $\rho_\infty\geq
0$ and free boundary $l_\infty>0$. If
    \begin{equation}
    \gamma=\frac{2n-2}{n}\ \textrm{ and }\ \ Gn^\frac{2-n}{n}M^\frac{2}{n}<2A
    \label{sym-self-E1.13-1}
    \end{equation}
or
     \begin{equation}
     \gamma>\frac{2n-2}{n},\label{sym-self-E1.14}
     \end{equation}
 from Proposition
      \ref{sym-self-stat-ex-prop5}, we know that there exists a unique
solution $(\rho_\infty,l_\infty)$ to the stationary system
(\ref{E-sym1})-(\ref{E-sym3}), satisfying $0<\underline{\rho}\leq
\rho_\infty(r)\leq \bar{\rho}<\infty$, $(\rho_\infty)_r(r)<0$,
$0<r<l_\infty$ with $l_\infty<+\infty$.

To handle the free boundary problem
(\ref{sym-E2.1})-(\ref{sym-E2.4}), it is convenient to reduce the
problem in Eulerian coordinates $(r,\tau)$ to the problem in
Lagrangian coordinates $(x,t)$, via the transformation:
    \begin{equation}
      x=\int^r_0y^{n-1}\rho(y,\tau)dy,\ t=\tau.\label{sym-self-E1.13}
    \end{equation}
Then the fixed boundary $r=0$ and the free boundary $r=b(\tau)$
become
        $$
          x=0\ \textrm{ and }x=\int^{b(\tau)}_0y^{n-1}\rho(y,\tau)dy
          =\int^{b}_0y^{n-1}\rho_0(y)dy=M,
        $$
where $M$ is the total mass initially. So that the region
$\{(r,\tau):0\leq r\leq b(\tau),\tau\geq0\}$ under consideration
is transformed into the region $\{(x,t):0\leq x\leq M, t\geq0\}$.

    Under the coordinate transformation (\ref{sym-self-E1.13}), the equations
    (\ref{sym-E2.1})-(\ref{sym-E2.4}) are transformed into
   \begin{equation}
      \left\{\begin{array}{l}
      \partial_t \rho(x,t)=-\rho^2\partial_x(r^{n-1} u),\\
           \partial_t u(x,t)=r^{n-1}\left\{\partial_x[\rho(\lambda+2\mu)\partial_{x}(r^{n-1}u)-P]
           -2(n-1)\frac{u}{r}\partial_x\mu\right\}-f(x,r,t),\\
      r^n(x,t)=n\int^x_0\rho^{-1}(y,t)dy,
     \end{array}
      \right.\label{sym-E}
    \end{equation}
 where $(x,t)\in(0,M)\times(0,\infty)$, with the initial data
    \begin{equation}
      (\rho,u)|_{t=0}=(\rho_0,u_0)(x),
        r|_{t=0}=r_0(x)=\left(
        n\int^x_0\rho^{-1}_0(y)dy
        \right)^\frac{1}{n},
    \end{equation}
and  the boundary conditions:
            \begin{equation}
            u(0,t)=0,
            \label{sym-Efixbd}
            \end{equation}
                \begin{equation}
                 \left.\left\{ P-\rho(\lambda+2\mu)\partial_x(r^{n-1}u)
                  +2(n-1)\mu\frac{u}{r}\right\}\right|_{x=M}=P_\Gamma,
                  \ t>0.
                  \label{sym-Efd}
                \end{equation}

                 It is standard that if we can solve the problem
(\ref{sym-E})-(\ref{sym-Efd}), then the free boundary problem
(\ref{sym-E1.1})-(\ref{sym-E1.4}) has a solution.

From (\ref{E-sym1})-(\ref{E-sym3}), it is easy to see that
$\rho_{\infty}(x)$ is the solution to the stationary system,
    \begin{equation}
      Ar_\infty^{n-1}(\rho_\infty^\gamma)_x=-f_{\infty}(x,r_{\infty}),
      \ r_{\infty}^n(x)=n\int^x_0\rho^{-1}_{\infty}(y)dy,
      \ x\in(0,M),\label{sym-E1.17}
    \end{equation}
        \begin{equation}
          \rho_\infty(M)=\left(\frac{P_\infty}{A}\right)^{\frac{1}{\gamma}}.\label{sym-E1.18}
        \end{equation}

The results in \cite{hoff,xin} show that the compressible
Navier-Stokes system with the constant viscosity coefficient have
the singularity at the vacuum.  Considering the modified
Navier-Stokes system in which the viscosity coefficient depends on
the density, Liu, Xin and Yang in \cite{liu} proved that such
system is local well-posedness. It is motivated by the physical
consideration that in the derivation of the Navier-Stokes
equations from the Boltzmann equation through the Chapman-Enskog
expansion to the second order, cf. \cite{Grad}, the viscosity
coefficient is a function of the temperature. If we consider the
case of isentropic fluids, this dependence is reduced to the
dependence on the density function.

 Since $n\geq2$ and the viscosity coefficient $\mu$ depends on $\rho$,
 the  nonlinear term  $2(n-1)\frac{1}{r}u\partial_x\mu$ in (\ref{sym-E})$_2$ makes the analysis
  significantly different from  the one-dimensional
case \cite{liu,okada1,vong,yang3,fang06}. Considering the
compressible spherically symmetric Navier-Stokes equations without
a solid core, the techniques in the case of similar system with a
solid core
\cite{Chen2002,Ducomet2005,MatusuNecasova,okada93,Zlotnik2005}
failed to be of use in our case, so we need obtain some new
\textit{a priori} estimates.

For spherically symmetric solutions of the Navier-Stokes equations
with constant viscosity, in \cite{hoff92}, the author gave an
information near the origin that the solution may develop vacuum
region about the origin. The difficulty of this problem is to
obtain the lower bound of the density $\rho$ and the upper bound
of the term  $\frac{1}{r}u$. When the initial data are small in
some sense, using some new \textit{a priori} estimates on the
solution, we can obtain the lower bound of the density and the
upper bound of the term  $\frac{1}{r}u$. The key ideas are using
the classical continuity method and the result of \textbf{Claim
1}. In \textbf{Claim 1}, we want to prove that there is a small
positive constant $\epsilon_1$, such that, for any $T>0$, if
    $$
      I(t)=\|\rho(\cdot,t)-\rho_\infty\|_{L^\infty}
      +\left\|\frac{u}{r}(\cdot,t)\right\|_{L^\infty}\leq
      2\epsilon_1,\ \forall\ t\in[0,T],
    $$
then
    $$
    I(t)\leq \epsilon_1, \ \forall\ t\in[0,T].
    $$
Let
    \begin{eqnarray*}
      B[\rho,u,r]&=&\int^M_0
      \left[(\rho-\rho_\infty)^2+r^{2n-2+\alpha}(\rho-\rho_\infty)_x^2+\frac{u^2}{r^2}
      \right.\nonumber\\
            &&+r^{2n-2}u_x^2+\left.
      r^{2n-2+\alpha} (\rho^{1+\theta}(r^{n-1}u)_x)_x^2\right]dx,
    \end{eqnarray*}
where $\alpha=\frac{3}{2}-n$. In Lemmas
\ref{sym-L2.2}-\ref{sym-L2.8}, we get some uniform \textit{a
priori} estimates (with respect to time) on the solution in the
weighted Sobolev space and the upper bound of $B[\rho,u,r]$. Using
the bound of  $B[\rho,u,r]$ and Sobolev's embedding Theorem, we
can finish the proof of \textbf{Claim 1}. Then, we will
 construct a weak solution  by using the finite
difference approximation. Our results show that: such system does
not develop vacuum states or concentration states for all time,
and the interface $\partial\Omega_\tau$ propagates with finite
speed. Since these estimates of the solution are uniform in time,
we could show that the solution converges to a stationary one as
time tends to infinity. Moreover, we construct various Lyapunov
functionals and obtained the stabilization rate estimates of
exponential type.

We now briefly review the previous works in this direction.
 For the related free boundary problem of one-dimensional isentropic
fluids with density-dependent viscosity (like
$\mu(\rho)=c\rho^\theta$), see \cite{liu,okada1,vong,yang3,fang06}
and the references therein. For the  spherically symmetric
solutions of the Navier-Stokes equations with
 a free boundary, see \cite{Chen2002,Ducomet2005,MatusuNecasova,okada93,Zlotnik2005}
 etc.. Ducomet-Zlotnik\cite{Ducomet2005,Zlotnik2005}
studied the similar system with a solid core and
 without the nonlinear term  $2(n-1)\frac{1}{r}u\partial_x\mu$.  Also see Lions\cite{Lions} and
Vaigant-Kazhikhov\cite{Vaigant} for multidimensional isentropic
fluids. For the related stabilization rate estimates in the
one-dimensional case, see
\cite{Ducomet2005-2,Matsumura,Mucha,Straskraba,fang06} etc..

Main assumptions on $c_1$, $c_2$, $\theta$ and $\gamma$ can be
stated as follows:
    \begin{description}
\item [(A1)]  condition (\ref{sym-self-E1.13-1}) or
(\ref{sym-self-E1.14}) holds;

\item [(A2)] $\theta\geq0$. $c_1$ and $c_2$ satisfy that
        $$
    c_1>0,\ 2c_1+nc_2>0
        $$
and
    \begin{equation}
    [2c_1\alpha+c_2(2n-2+\alpha)]^2-4(2c_1+c_2)[2c_1(n-1)+c_2(n-1)(n-1+\alpha)]<0,\label{sym-E1.19}
    \end{equation}
where $\alpha=\frac{3}{2}-n$.
    \end{description}

Under the above assumptions (A1)-(A2), we will prove the existence
of a global weak solution to the initial-boundary value problem
(\ref{sym-E})-(\ref{sym-Efd}) in the sense of the following
definition.
\begin{defn} A pair of functions $(\rho,u,r)(x,t)$ is called a
global weak solution to the initial boundary value problem
(\ref{sym-E})-(\ref{sym-Efd}), if for any $T>0$,
    $$
    \rho,u \in L^\infty([0,M]\times[0,T])\cap C^1([0,T];L^2([0,M])),
    $$
        $$
          r\in C^{1}([0,T];L^\infty([0,M])),
        $$
        $$
       (r^{n-2}u)_x,   (r^{n-1})_x\in L^\infty([0,T];L^{n-\frac{1}{2}}([0,M])),
        $$
 and
            $$
         (r^{n-1}u)_x\in L^\infty([0,M]\times[0,T])\cap C^\frac{1}{2}([0,T];L^2([0,M])).
        $$
Furthermore, the following equations hold:
    $$
    \rho_t+\rho^2(r^{n-1}u)_x=0,\
    \rho(x,0)=\rho_0(x)\ a.e.
    $$
    $$
      r_t=u,\ r^n(x,t)=n\int^x_0\rho^{-1}(y,t)dy,\ r(x,0)=r_0(x)\ a.e.
    $$
and
        \begin{eqnarray*}
    &&\int^\infty_0\int^M_0[u\psi_t+(P-\rho(\lambda+2\mu)(r^{n-1}
    u)_x)(r^{n-1}\psi)_x\nonumber\\
    &&+2(n-1)\mu(r^{n-2}u\psi)_x-f(x,r,t)\psi]dxdt\nonumber\\
        &=&\int^\infty_0P_\Gamma(r^{n-1}\psi)(M,t)dt-\int^M_0u_0(x)\psi(x,0)dx,
        \end{eqnarray*}
for any test function $\psi(x,t)\in C^\infty_0(\Omega)$ with
$\Omega=\{(x,t):\ 0< x\leq M,\ t\geq 0\}$.
\end{defn}

In what follows, we always use $C$($C_i$) to denote a generic
positive constant depending only on the initial data, independent
of the given time T.

We now state the main theorems in this paper. Let
$\underline{\rho}=\min_{x\in[0,M]}\rho_\infty$ and
$\bar{\rho}=\max_{x\in[0,M]}\rho_\infty$.
\begin{thm}\label{sym-thm} Under the conditions (\ref{sym-f1})-(\ref{f4}) and [A1]-[A2], there exists
a positive constant $\epsilon_0>0$, such that if
    \begin{equation}
             \|(f_1,\Delta P)\|_{L^\infty\cap
      L^1}+\|(\Delta P)'\|_{L^2}+\|f_2\|_{L^2}\leq \epsilon_0,\label{epsilon1}
    \end{equation}
    \begin{equation}
      \|\rho_0-\rho_\infty\|_{L^\infty}^2+ B[\rho_0,u_0,r_0]\leq
      \epsilon_0^2,\label{epsilon2}
    \end{equation}
 then the system
(\ref{sym-E})-(\ref{sym-Efd}) has a unique global weak solution
$(\rho,u,r)$ satisfying
 \begin{equation}
          \rho(x,t)\in \left[\frac{1}{2}\underline{\rho},
          \frac{3}{2}\bar{\rho}
          \right],\ r^n(x,t)\in [C^{-1}x,Cx],\label{sym-E1.32}
        \end{equation}
  \begin{equation}
   \left\|\frac{u}{r}(\cdot,t)
   \right\|_{L^\infty}\leq C \|\partial_x(r^{n-1}u)(\cdot,t)\|_{L^\infty}
  \leq C\epsilon_0,\label{sym-E1.33-1}
  \end{equation}
    \begin{equation}
      B[\rho,u,r]\leq C\epsilon_0^2,\label{sym-E1.40-1}
    \end{equation}
for all $t\geq0$ and $x\in[0,M]$. Furthermore, we have
     $$
    \lim_{t\rightarrow+\infty}\int^M_0\left\{x^{\frac{2n-2+\alpha}{n}}u_x^2+
    x^{\frac{2n-2+\alpha}{n}}\left[(\rho^\theta)_x-(\rho_\infty^\theta)_x\right]^2
    \right\}dx=0,
        $$
        $$
    \lim_{t\rightarrow+\infty}\|u(\cdot,t)\|_{L^\infty}+
    \|\rho(\cdot,t)-\rho_\infty(\cdot)\|_{L^\infty}+
        \|r(\cdot,t)-r_\infty(\cdot)\|_{L^\infty}=0.
        $$
\end{thm}

\begin{rem}
In fact, assumption (\ref{sym-E1.19}) is give a restriction on
$\frac{\lambda}{\mu}$, i.e.
    \begin{eqnarray*}
        &&\frac{-18+8n+8n^2-8\sqrt{3}(n-1)\sqrt{4n-3}}{9-12n+4n^2}\\
        &<&\frac{\lambda}{\mu}
        <
        \frac{-18+8n+8n^2+8\sqrt{3}(n-1)\sqrt{4n-3}}{9-12n+4n^2}.
    \end{eqnarray*}
If $n=3$, we can choose
$\frac{\lambda}{\mu}=\frac{c_2}{c_1}\in(\frac{2}{3}(13-8\sqrt{3}),\frac{2}{3}(13+8\sqrt{3}))$.
\end{rem}
\begin{rem}
We can choose  the constant $\epsilon_0$ as in
(\ref{sym-self-E3.77}).
\end{rem}

The proof of the uniqueness part  of Theorem \ref{sym-thm} also
shows that the continuous dependence of the solution on the
initial data holds. We may state the following result without a
proof.

\begin{thm}
For each $i=1,2$, let $(\rho_i,u_i,r_i)$ be the solution to the
system (\ref{sym-E})-(\ref{sym-Efd}) with the initial data
$(\rho_{0i},u_{0i},r_{0i})$, which satisfy regularity conditions
(\ref{sym-E1.32})-(\ref{sym-E1.40-1}). Then, we have
    \begin{eqnarray*}
    &&\int^M_0[(u_1-u_2)^2+(\rho_1-\rho_2)^2+x^{-\frac{2}{n}}(r_1-r_2)^2](x,t)dx\nonumber\\
        &\leq& Ce^{Ct}
        \int^M_0[(u_{01}-u_{02})^2+(\rho_{01}-\rho_{02})^2+x^{-\frac{2}{n}}(r_{01}-r_{02})^2]dx,
    \end{eqnarray*}
for all $t\geq0$.
\end{thm}

\begin{thm}\label{sym-thm2}
  Under the assumptions of Theorem \ref{sym-thm} and
    \begin{equation}
      f_1(t)
      +f_2(t)+|\Delta P(t)|+|(\Delta P)'(t)|\leq Ce^{-a_0t},\label{sym-E1.42}
    \end{equation}
where $a_0$ is a positive constant, then we have
    $$
      \int^M_0\left\{r^{2n-2+\alpha}(\rho-\rho_\infty)_x^2
      +r^{2n-2+\alpha}[\partial_x(\rho^{1+\theta}\partial_x(r^{n-1}u))]^2+r^\alpha u_t^2\right\}dx\leq
      Ce^{-at},
    $$
    $$
  \left\|\left(\frac{u}{r},
     (r^{n-1}u)_x\right)(\cdot,t)\right\|_{L^\infty}+ \|\rho(\cdot,t)-\rho_\infty(\cdot)\|_{L^\infty}
      +\|r(\cdot,t)-r_\infty(\cdot)\|_{L^\infty}\leq Ce^{-at},
    $$
for all $t\geq0$, where $a$ is a positive constant.
\end{thm}

\begin{rem}
Considering the general case that $(\mu,\lambda)(\cdot)\in
C(\mathbb{R_+})\cap W^{1,\infty}_{loc}(\mathbb{R_+})$, under the
assumptions (\ref{sym-f1})-(\ref{f4}), (A1) and
    $$
    \mu(\rho)>0,\ 2\mu(\rho)+n\lambda(\rho)>0,
        $$
    $$
      [2\mu\alpha+\lambda(2n-2+\alpha)]^2-4(2\mu+\lambda)
      [2\mu(n-1)+\lambda(n-1)(n-1+\alpha)]<0,
    $$
 for all
$\rho\in[\frac{1}{2}\underline{\rho},\frac{3}{2}\bar{\rho}]$, we
can obtain  the same results.
\end{rem}
\begin{rem}
In this paper, we study the case of $\gamma>1$ and prove the main
results in this case only, since the case of $\gamma=1$ can be
discussed through the similar process. The main different is that
(\ref{sym-self-E2.9-1}) is replaced by
    \begin{equation*}
        S[V]=\int^M_0\left(A\ln V_x+P_\infty V_x
        +\int^V_1Gx(nh)^{\frac{2-2n}{n}}dh
    \right)dx,
    \end{equation*}
when $\gamma=1$ and $n=2$.
\end{rem}

The rest of this paper is organized as follows. First, we obtain
the existence and uniqueness of the solution to the stationary
problem in Section \ref{sym-self-sec1.5}. In Section
\ref{sym-sec2}, we will prove some \textit{a priori} estimates
which will be used to obtain global existence of the weak
solutions. In Section \ref{sym-Sec3}, using the finite difference
approximation and  \textit{a priori} estimates obtained in Section
\ref{sym-sec2}, we prove the existence part of Theorem
\ref{sym-thm}. In Section \ref{sym-Sec4}, we will prove the
uniqueness of the weak solution. In Section \ref{sym-Sec5}, We
 show that the solution of the free boundary problem tends to
a stationary one, as $t\rightarrow+\infty$. In Section
\ref{sym-Sec6}, we will obtain the stabilization rate estimates of
exponential type on the solution by constructing some Lyapunov
functionals.

\section{The stationary problem}\label{sym-self-sec1.5}
We start with a proof of the existence of a positive solution to
the Lagrangian stationary problem.
Zlotnik-Ducomet\cite{Zlotnik2005} studied the stationary problem
with a solid core $r\geq r_0>0$. Using similar arguments as that
in \cite{Zlotnik2005}, we can obtain the following results for the
stationary problem without a solid core.

\begin{prop}\label{sym-self-stat-ex-prop1}
If
    \begin{equation}
    \gamma>\frac{2n-2}{n}\label{sym-self-E2.1}
    \end{equation}
or
    \begin{equation}
      \gamma=\frac{2n-2}{n}
      \ \textrm{ and }\ Gn^\frac{2-n}{n}M^\frac{2}{n}<2A,\label{sym-self-E2.2}
    \end{equation}
or
    \begin{equation}
    0<\gamma<\frac{2n-2}{n}
      \ \textrm{ and }\
      P_\infty+\frac{G}{2}n^\frac{2-n}{n}M^\frac{2}{n}\delta_3^\frac{2n-2}{n}\leq A\delta_3^\gamma,
    \end{equation}
where $\delta_3=\left(
      \frac{A\gamma n^\frac{2n-2}{n}}{(n-1)GM^{\frac{2}{n}}}
      \right)^\frac{n}{2n-2-n\gamma}$,
 then the Lagrangian stationary problem
(\ref{sym-E1.17})-(\ref{sym-E1.18}) has a positive solution
$\rho_\infty\in W^{1,\beta}([0,M])$, where
$\beta\in[1,\frac{n}{n-2})$ is a constant.
\end{prop}
\begin{proof}
  We introduce the nonlinear operator
    $$
    I:K\rightarrow W^{1,\beta}([0,M]),
    $$
where
    $K=\left\{f\in C([0,M]):\displaystyle{\min_{x\in[0,M]}}f(x)\geq
    \left(\frac{P_\infty}{A}\right)^\frac{1}{\gamma}
    \right\}$, by setting
        $$
        I(f)(x)=\left(\frac{P_\infty+\int^M_xG\frac{y}{r_f^{2n-2}(y)}dy}{A}
        \right)^\frac{1}{\gamma}
        $$
with $r_f^n(x)=n\int^x_0f^{-1}(y)dy$,
    $x\in[0,M]$.
We can restate the problem (\ref{sym-E1.17})-(\ref{sym-E1.18}) as
the fixed-point problem
    \begin{equation}
    \rho_\infty=I (\rho_\infty).\label{sym-self-E2.3}
    \end{equation}

For all $f\in K_\delta=\left\{ f\in K: f\leq \delta\right\}$ with
$\delta>\left(\frac{P_\infty}{A}\right)^\frac{1}{\gamma}$, we have
        $$
        nx\delta^{-1}\leq r_f^n(x)
        $$
and
    \begin{eqnarray*}
      P_\infty\leq A(I(f))^\gamma&\leq& P_\infty
      +G\delta^{\frac{2n-2}{n}}n^{-\frac{2n-2}{n}}\int^M_0x^{\frac{2-n}{n}}dx\nonumber\\
        &=&P_\infty
      +\frac{G}{2}\delta^{\frac{2n-2}{n}}n^{\frac{2-n}{n}}M^\frac{2}{n}.
    \end{eqnarray*}

If $\gamma>\frac{2n-2}{n}$, then $I (K_{\delta_1})\subset
K_{\delta_1}$, where $\delta_1$ is a positive constant satisfying
    $P_\infty+\frac{G}{2}\delta_1^{\frac{2n-2}{n}}n^{\frac{2-n}{n}}M^\frac{2}{n}\leq
    A\delta_1^\gamma$. And one can immediately verify that $I$ is a
    compact operator on $K_{\delta_1}$. Since $K_{\delta_1}$  is a convex
 closed bounded non-empty subset of $C([0,M])$, the problem
 (\ref{sym-self-E2.3}) has a solution $\rho\in K_{\delta_1}$ by Schauder's fixed point theorem.

If $\gamma=\frac{2n-2}{n}$ and $Gn^\frac{2-n}{n}M^\frac{2}{n}<2A$,
then $I (K_{\delta_2})\subset K_{\delta_2}$, where $\delta_2$ is a
positive constant satisfying
    $P_\infty+\frac{G}{2}\delta_2^{\frac{2n-2}{n}}n^{\frac{2-n}{n}}M^\frac{2}{n}\leq
    A\delta_2^\gamma$.

If $\gamma<\frac{2n-2}{n}$ and
 \begin{equation*}
          P_\infty+\frac{G}{2}n^\frac{2-n}{n}M^\frac{2}{n}\left(
      \frac{A\gamma n^\frac{2n-2}{n}}{(n-1)GM^{\frac{2}{n}}}
      \right)^\frac{2n-2}{2n-2-n\gamma}\leq A\left(
      \frac{A\gamma n^\frac{2n-2}{n}}{(n-1)GM^{\frac{2}{n}}}
      \right)^\frac{n\gamma}{2n-2-n\gamma}
    \end{equation*}
then $I( K_{\delta_3})\subset K_{\delta_3}$, where
    $$\delta_3=\left(
      \frac{A\gamma n^\frac{2n-2}{n}}{(n-1)GM^{\frac{2}{n}}}
      \right)^\frac{n}{2n-2-n\gamma}.$$
We can finish the proof of the theorem immediately.
\end{proof}

Letting $V_\infty=\frac{r^n_\infty}{n}$, using the equality
$\frac{1}{\rho_\infty}=(V_\infty)_x$, one can eliminate the
function $\rho_\infty$ from the Lagrangian stationary problem
(\ref{sym-E1.17})-(\ref{sym-E1.18}) and obtain an equivalent
boundary-value problem for a non-linear second-order ODE:
    \begin{equation}
      (A(V_\infty)_x^{-\gamma})_x=-Gxn^{\frac{2-2n}{n}}V_\infty^{\frac{2-2n}{n}},\
      x\in(0,M),\label{sym-self-E2.4}
    \end{equation}
        \begin{equation}
          V_\infty(0)=0,\
          (V_\infty)_x(M)=\left(\frac{A}{P_\infty}\right)^{\frac{1}{\gamma}},
          \label{sym-self-E2.5}
        \end{equation}
for a function $V_\infty\in C^{1}([0,M])$ such that
$(V_\infty)_x>0$.

    In accordance with the method of small perturbations, we
    replace $V_\infty$ by $V=V_\infty+W$ with small $W$ and
    linearized the operator in the last problem:
    \begin{eqnarray*}
    &&(A(V)_x^{-\gamma})_x+Gxn^{\frac{2-2n}{n}}V^{\frac{2-2n}{n}}\\
        &=&
(-\gamma A(V_\infty)_x^{-\gamma-1}W_x)_x+(2-2n)Gx(n
V_\infty)^{\frac{2-3n}{n}}W+\ldots, \ x\in(0,M),
    \end{eqnarray*}
        $$
    V(0)=0+W(0),\ A(V_x)^{-\gamma}|_{x=M}-P_\infty=-\gamma
    A\{(V_\infty)_x^{-\gamma-1}W_x\}|_{x=M}+\ldots,
        $$
up to the terms of the second order of smallness with respect to
$W$. We define the linearized operator
    \begin{equation}
    L[W]=(-\gamma A\rho_\infty^{\gamma+1}W_x)_x+(2-2n)Gx(n
V_\infty)^{\frac{2-3n}{n}}W,\ W\in K_0,\label{sym-self-E2.6}
    \end{equation}
where
        $
     K_0=\{W\in C^1([0,M]):     W(0)=0, W_x(M)=0\}.
        $
It is easy to get
    \begin{equation*}
    (L[W],W)=
      \int^M_0\left(\gamma A(\rho_\infty)^{1+\gamma}W_x^2
    -(2n-2)Gx(nV_\infty)^{\frac{2-3n}{n}} W^2
    \right)dx,\ W\in K_0.
    \end{equation*}
Let
    \begin{equation}
    J[W]:=
      \int^M_0\left(\gamma A(\rho_\infty)^{1+\gamma}W_x^2
    -(2n-2)Gx(nV_\infty)^{\frac{2-3n}{n}} W^2
    \right)dx,
    \end{equation}
for $W\in K_1=\{f\in C^1([0,M]): f(0)=0\}$.

 We say a stationary solution $V_\infty$ is \textit{statically
stable} if
    \begin{equation}
      J[W]\geq \delta_3\left(\|W_x(x)\|_{L^2(0,M)}^2+
      \|x^{-1}W(x)\|_{L^2(0,M)}^2
      \right),\label{sym-selfE2.8}
    \end{equation}
for some $\delta_3>0$ and all $W\in K_1$.

Now, the static potential energy  takes the following form:
    \begin{equation}
        S[V]=\int^M_0\left(\frac{A}{\gamma-1}(V_x)^{1-\gamma}+P_\infty V_x
        +\int^V_1Gx(nh)^{\frac{2-2n}{n}}dh\label{sym-self-E2.9-1}
    \right)dx.
    \end{equation}
We call  $V\in K_2=\{f\in C^1([0,M]):f(0)=0, \min(f_x)>0\}$ is a
point of \textit{local quadratic minimum of $S$} if
    \begin{equation}
      S[V+W]-S[V]\geq \delta_4\left(\|W_x(x)\|_{L^2(0,M)}^2+
      \|x^{-1}W(x)\|_{L^2(0,M)}^2
      \right),\label{sym-self-E2.9}
    \end{equation}
for all $W\in K_1$ and $\|W\|_{C^1([0,M])}\leq \delta_5$, for some
$\delta_4>0$ and $\delta_5>0$.

We can clarify the variational sense of the definition of
\textit{statically stable} as follows.
\begin{prop}\label{sym-self-stat-ex-prop2}
A function $V\in K_2$ is a point of local quadratic minimum of $S$
if and only if $V=V_\infty$ is a solution of the problem
(\ref{sym-self-E2.4})-(\ref{sym-self-E2.5}) and satisfies static
stability condition (\ref{sym-selfE2.8}).
\end{prop}
\begin{proof}
  Let $V\in K_2$, $W\in K_1$ and $\|W\|_{C^1([0,M])}=1$. Using Taylor's formula,
  we have
        $$
    S[V+\epsilon W]=S[V]+\delta
    S[V](\epsilon W)+\frac{1}{2}\frac{d^2}{d\tau^2}S[V+\tau\epsilon
    W]\big|_{\tau=\tilde{\tau}},
        $$
  where
    $$
    \delta S[V](\epsilon W)=\int^M_0\left(-A(V_x)^{-\gamma}\epsilon W_x
    +P_\infty \epsilon W_x+Gx(nV)^{\frac{2-2n}{n}}\epsilon W
    \right)dx,
    $$
  and
    \begin{eqnarray*}
    \frac{d^2}{d\tau^2}S[V+\tau\epsilon
    W]&=&\int^M_0\left(\gamma A(V_x+\tau\epsilon W_x)^{-1-\gamma}(\epsilon W_x)^2
    \right.\\
        &&\left. -(2n-2)Gx(n(V+\tau\epsilon W))^{\frac{2-3n}{n}}(\epsilon W)^2
    \right)dx,
    \end{eqnarray*}
for all  $|\epsilon|<\frac{1}{\min V_x}$ and some
$\tilde{\tau}\in[0,1]$. If (\ref{sym-self-E2.9}) holds, we have
    $$
    \frac{d^2}{d\tau^2}S[V+\tau\epsilon
    W]\leq C\epsilon^2\left(\|W_x(x)\|_{L^2(0,M)}^2+
      \|x^{-1}W(x)\|_{L^2(0,M)}^2
      \right)
    $$
and
    $$
    C\epsilon^2\left(\|W_x(x)\|_{L^2(0,M)}^2+
      \|x^{-1}W(x)\|_{L^2(0,M)}^2
      \right)+\epsilon \delta S[V](W)>0,
    $$
 for all $|\epsilon|\in(0,\min(\delta_5,\frac{1}{\min
    V_x}))$ and $\|W\|_{C^1([0,M])}=1$.
Thus, we obtain
    $$
    \delta S[V](W)=0,
    $$
i.e.
    $$
    \int^M_0\left(-A(V_x)^{-\gamma} W_x
    +P_\infty  W_x+Gx(nV)^{\frac{2-2n}{n}}W
    \right)dx=0,
    $$
 for all $W\in K_1$ and $\|W\|_{C^1{[0,M]}}=1$, that is, $V$ is a stationary
 point of $S$ and a solution of the problem
 (\ref{sym-self-E2.4})-(\ref{sym-self-E2.5}).
We can rewrite $\frac{d^2}{d\tau^2}S[V+\tau\epsilon
    W]$ as follows
    $$
    \frac{d^2}{d\tau^2}S[V+\tau\epsilon
    W]=\delta^2S[V](\epsilon W)+S_1,
    $$
where $\delta^2S[V](\epsilon
W)=\frac{d^2}{d\tau^2}S[V+\tau\epsilon
    W]\big|_{\tau=0}$ and
        \begin{eqnarray*}
         | S_1|&=&\left|\frac{d^2}{d\tau^2}S[V+\tau\epsilon
    W]-\delta^2S[V](\epsilon W)\right|\\
            &\leq& C\epsilon\left(\|\epsilon
            W_x(x)\|^2_{L^2([0,M])}+\|x^{-1}\epsilon W(x)\|^2_{L^2([0,M])}
            \right).
        \end{eqnarray*}
Thus, we obtain
    $$
    \delta^2S[V](\epsilon W)\geq (\delta_4-C\epsilon) \left(\|\epsilon
            W_x(x)\|^2_{L^2([0,M])}+\|x^{-1}\epsilon W(x)\|^2_{L^2([0,M])}
            \right)
    $$
for all $\epsilon\in(0,\min(\delta_5,\frac{1}{\min
    V_x},\frac{\delta_4}{2C}))$ and $\|W\|_{C^1([0,M])}=1$.
Moreover, we have
     \begin{equation}
   J[W]:=\delta^2S[V](W)\geq \frac{\delta_4}{2}\left(\|
            W_x(x)\|^2_{L^2([0,M])}+\|x^{-1}  W(x)\|^2_{L^2([0,M])}
            \right),
    \end{equation}
for all $W\in K_1$.

 If $V=V_\infty$ is a solution of
the problem (\ref{sym-self-E2.4})-(\ref{sym-self-E2.5}) and
satisfies static stability condition (\ref{sym-selfE2.8}), we can
prove $V_\infty$ is a point of local quadratic minimum of $P$
easily.
\end{proof}

\begin{prop}\label{sym-self-stat-ex-prop3}
If $V=V_\infty$ is a solution of the problem
(\ref{sym-self-E2.4})-(\ref{sym-self-E2.5}) and
$\gamma\geq\frac{2n-2}{n}$, then
  (\ref{sym-selfE2.8}) and (\ref{sym-self-E2.9}) hold.
\end{prop}
\begin{proof}
  From $(A\rho_\infty^\gamma)_x=-G\frac{x}{r_\infty^{2n-2}}=-Gx(nV_\infty)^{\frac{2-2n}{n}}$,
  using integration by parts,  we
  have
    \begin{eqnarray}
    J[W]  &=&\int^M_0\left(\gamma A(\rho_\infty)^{1+\gamma}W_x^2
    -(2n-2)Gx(nV_\infty)^{\frac{2-3n}{n}} W^2
    \right)dx\nonumber\\
            &=&\int^M_0\left(\gamma A(\rho_\infty)^{1+\gamma}W_x^2
    +(2n-2)A(\rho_\infty^\gamma)_x(nV_\infty)^{-1} W^2
    \right)dx\nonumber\\
            &=&\int^M_0(\gamma A(\rho_\infty)^{1+\gamma}W_x^2
    -2(2n-2)A\rho_\infty^\gamma(nV_\infty)^{-1} WW_x
    \nonumber\\
    &&+\frac{2n-2}{n}A\rho_\infty^{\gamma-1}V_\infty^{-2} W^2
    )dx+(2n-2)P_\infty\left(\frac{W^2}{nV_\infty}\right)(M)\nonumber\\
    &:=&I_0[W]+(2n-2)P_\infty\left(\frac{W^2}{nV_\infty}\right)(M),
    \ \textrm{for all }\ W\in K_0.\label{sym-self-E2.12}
    \end{eqnarray}
If $\gamma\geq\frac{2n-2}{n}$, we have
    \begin{equation}
      I_0[W]\geq
      \int^M_0\frac{2n-2}{n}A\rho_\infty^{1+\gamma}\left(W_x-\frac{W}{\rho_\infty V_\infty}\right)^2dx.
    \end{equation}
If (\ref{sym-selfE2.8}) not holds, we have for any integer $m>1$,
there exists $W_m\in K_0$ and $\|W_m\|_{C^1([0,M])}=1$ such that
        \begin{equation}
    J[W_m]<\frac{1}{m}\left(\|(W_m)_x(x)\|_{L^2(0,M)}^2+
      \|x^{-1}W_m(x)\|_{L^2(0,M)}^2
      \right).      \label{sym-self-E2.13}
        \end{equation}
Then, there is a subsequence $m\rightarrow\infty$ for which
    $$
    W_m\rightarrow W\ \textrm{ in }\ C([0,M]),
    $$
        $$
        (W_{m})_x\rightharpoonup W_x \ \textrm{ in }\ L^2([0,M]).
        $$
From (\ref{sym-self-E2.12})-(\ref{sym-self-E2.13}), we have
    $$
    W_x=\frac{W}{\rho_\infty V_\infty},\ x\in(0,M),
    $$
 and $W(0)=W(M)=0$. Thus, we
obtain $W\equiv0$. It is a contradiction.

Therefore, if $\gamma\geq\frac{2n-2}{n}$, then
(\ref{sym-selfE2.8}) holds. From Proposition
\ref{sym-self-stat-ex-prop1}-\ref{sym-self-stat-ex-prop2}, we can
obtain (\ref{sym-self-E2.9}) immediately.
\end{proof}

Now, we shall use the shooting method to prove the uniqueness of
the solution.
\begin{prop}\label{sym-self-stat-ex-prop4}
Under the assumptions (\ref{sym-self-E2.1})-(\ref{sym-self-E2.2}),
the Lagrangian stationary problem
(\ref{sym-E1.17})-(\ref{sym-E1.18}) has a unique positive solution
$\rho_\infty$.
\end{prop}
\begin{proof}
 We
consider the Cauchy problem
    \begin{equation}
      (A\rho^\gamma_\infty)_x=-Gx(nV_\infty)^{\frac{2-2n}{n}},
      \ (V_\infty)_x=\rho_\infty^{-1},
      \ x\in(0,M),\label{sym-self-E2.14}
    \end{equation}
        \begin{equation}
          \rho_\infty\big|_{x=0}=\sigma,
          \ V_\infty\big|_{x=0}=0,\label{sym-self-E2.15}
        \end{equation}
for the unknown functions $\rho_\infty(\sigma,x)$ and
$V_\infty(\sigma,x)$, where $\sigma>0$ is the shooting parameter.
For each $\sigma>0$, using  similar arguments as that in
Proposition \ref{sym-self-stat-ex-prop1}, we can obtain the
existence of the solution to this problem, satisfying
    \begin{equation}
    \rho_\infty(\sigma,x)\in [\frac{\sigma}{2},\sigma],
    \ V_\infty(\sigma,x)\in[\frac{x}{\sigma},\frac{2x}{\sigma}],
    \  x\in[0,M_0],\label{sym-self-E2.16}
    \end{equation}
    \begin{equation}
    \rho_\infty\in W^{1,\beta}([0,M_0]), V_\infty\in C^1([0,M_0]),\label{sym-self-E2.17}
    \end{equation}
where $M_0$ is a positive constant satisfying
    $A\sigma^\gamma-\sigma^{\frac{2n-2}{n}}\frac{G}{2}n^\frac{2-n}{n}M_0^\frac{2}{n}\geq A
    \left(\frac{\sigma}{2}\right)^\gamma$ and $M_0\leq M$.
If there exist two solutions $(\rho_1,V_1)$ and $(\rho_2,V_2)$ to
this problem satisfying
    \begin{equation}
      \rho_i\in W^{1,\beta}([0,M_i]),
      \ x\in [0,M_i],\label{sym-self-E2.20-1}
    \end{equation}
where $M_i\in(0,M]$, $i=1,2$. From (\ref{sym-self-E2.20-1}), there
exists a positive constant $M_3\in (0,\min\{M_1,M_2\})$ such that
    $$
    \rho_i(x)\in [\frac{\sigma}{2},\sigma]
    \ \textrm{ and } V_i(x)\in [\frac{x}{\sigma},\frac{2x}{\sigma}],
    \ x\in [0,M_3],\ i=1,2.
    $$
Then,  we have
    $$
      A\rho_1^\gamma-A\rho_2^\gamma=\int^x_0Gyn^\frac{2-2n}{n}(V_2^{\frac{2-2n}{n}}-V_1^{\frac{2-2n}{n}})dy
            \leq                   C\int^x_0y^{\frac{2-2n}{n}}\int^y_0|\rho_1^{-1}-\rho_2^{-1}|(z)dzdy,
    $$
and
    $$
     \|\rho_1-\rho_2\|_{L^\infty([0,\epsilon])}\leq C\|\rho_1-\rho_2\|_{L^\infty([0,\epsilon])}
     \int^\epsilon_0y^{\frac{2-n}{n}}dy\leq
     C_\sigma\epsilon^\frac{2}{n}\|\rho_1-\rho_2\|_{L^\infty([0,\epsilon])},
    $$
for all  $x,\epsilon\in(0,M_3]$. Choosing
$\epsilon<C_\sigma^{-\frac{n}{2}}$, we have
    $$
    \rho_1=\rho_2, \ \textrm{ for all }\ x\in[0,\epsilon].
    $$
Considering the Cauchy problem
    \begin{equation}
      (A\rho^\gamma_\infty)_x=-Gx(nV_\infty)^{\frac{2-2n}{n}},
      \ (V_\infty)_x=\rho_\infty^{-1},
      \ x\in(\frac{\epsilon}{2},M),
    \end{equation}
        \begin{equation}
          \rho_\infty\big|_{x=\frac{\epsilon}{2}}=\rho_1(\sigma,\frac{\epsilon}{2}),
          \ V_\infty\big|_{x=\frac{\epsilon}{2}}=\int^\frac{\epsilon}{2}_0\rho_1^{-1}(\sigma,y)dy,
        \end{equation}
using  the classical ODE theory, we have $\rho_1(x)=\rho_2(x)$,
 $x\in[\frac{\epsilon}{2},\min\{M_1,M_2\}]$. Thus, for each
$\sigma>0$, there exists a unique solution to the problem
(\ref{sym-self-E2.14})-(\ref{sym-self-E2.15}) satisfying
$\rho_\infty(x,\sigma)>0$ for $x\in[0,M_\sigma)$, where either
$\rho_\infty\big|_{x=M_\sigma}=0$ and $M_\sigma\in(0,M)$ or
$M_\sigma=M$.

Clearly, if $\rho_\infty$ is a solution to the problem
(\ref{sym-E1.17})-(\ref{sym-E1.18}), then $\rho_\infty$ satisfying
(\ref{sym-self-E2.14})-(\ref{sym-self-E2.15}) for some $\sigma>0$.
We will show that this can be possible only for one value of
$\sigma$. Using  similar arguments as that in the above part and
in \cite{Hartman} (\S V.3), we obtain that $(\partial_\sigma
\rho_\infty^\gamma,\partial_\sigma V_\infty)$ is  well defined and
satisfies the linear Cauchy problem
    \begin{equation}
      A(\partial_\sigma
\rho_\infty^\gamma)_x=(2n-2)Gx(nV_\infty)^\frac{2-3n}{n}\partial_\sigma
V_\infty, \ (\partial_\sigma
V_\infty)_x=-\frac{1}{\gamma}\rho_\infty^{-\gamma-1}\partial_\sigma
\rho_\infty^\gamma, \label{sym-self-E2.20}
    \end{equation}
where $ x\in[0,M_\sigma),$        \begin{equation}
    \partial_\sigma
\rho_\infty^\gamma\big|_{x=0}=1,\ \partial_\sigma
V_\infty\big|_{x=0}=0.\label{sym-self-E2.21}
        \end{equation}
It is easy to see that
    $$
    \partial_\sigma
\rho_\infty^\gamma>0,\ (\partial_\sigma V_\infty)_x<0,
\partial_\sigma V_\infty<0
        $$
hold on $[0,M_4)$, where either $\partial_\sigma
\rho_\infty^\gamma\big|_{x=M_4}=0$ and $M_4\in(0,M_\sigma)$ or
$M_4=M_\sigma$. We claim that only $M_4=M_\sigma$ can occur.

Assume that $M_4\in(0,M_\sigma)$. Letting
$\phi=A\rho_\infty^\gamma(\partial_\sigma
V_\infty)_x+\frac{n}{2n-2}A\partial_\sigma
\rho_\infty^\gamma(V_\infty)_x$, from (\ref{sym-self-E2.14}) and
(\ref{sym-self-E2.20}), we have
    $$
    \int^{M_4}_0\phi dx=\left.\left\{A\rho_\infty^\gamma\partial_\sigma
V_\infty+\frac{n}{2n-2}A\partial_\sigma \rho_\infty^\gamma
V_\infty\right\}\right|_0^{M_4}.
    $$
By the estimates $\rho_\infty(\sigma,M_4)>0$, $\partial_\sigma
\rho_\infty^\gamma\big|_{x=M_4}=0$, $\partial_\sigma
V_\infty\big|_{x=M_4}<0$ and the initial condition
(\ref{sym-self-E2.15}) and (\ref{sym-self-E2.21}), we get
    $$
    \int^{M_4}_0\phi dx<0.
    $$
On the other hand, from (\ref{sym-self-E2.14}) and
(\ref{sym-self-E2.20}), we have
    $$
    \phi=A\rho_\infty^{-1}\partial_\sigma\rho_\infty^\gamma(\frac{n}{2n-2}-\frac{1}{\gamma})\geq0,
    \ x\in (0,M_4).
    $$
It is a contradiction.

Thus, we obtain
    $$
    \rho_\infty>0,\ \partial_\sigma \rho_\infty>0,
    \ x\in (0,M_\sigma).
    $$
and $M_\sigma$ is non-decreasing on $\sigma\in (0,\infty)$.
Therefore, for each fixed point $x\in
[0,\sup_{\sigma>0}M_\sigma)$, the function $\rho_\infty(\sigma,x)$
is strictly increasing on $\sigma>\left(
\frac{P_\infty}{A}\right)^{\frac{1}{\gamma}}$, and satisfying
$A\rho_\infty^\gamma\big|_{x=M}=P_\infty$ for at most one value of
$\sigma$.
\end{proof}

Using the properties of the transformation (\ref{sym-self-E1.13})
and Propositions
\ref{sym-self-stat-ex-prop1}-\ref{sym-self-stat-ex-prop4}, we can
obtain the following proposition immediately.
\begin{prop}\label{sym-self-stat-ex-prop5}
  Under the assumptions (\ref{sym-self-E2.1})-(\ref{sym-self-E2.2}),
the Eulerian stationary problem (\ref{E-sym1})-(\ref{E-sym3}) has
a unique positive solution $(\rho_\infty,l_\infty)$, satisfying
$0<\underline{\rho}\leq \rho_\infty(r)\leq \bar{\rho}<\infty$,
$(\rho_\infty)_r(r)<0$, $0<r<l_\infty$ with $l_\infty<+\infty$.
\end{prop}
\section{\textit{A priori} estimates}\label{sym-sec2}

From (\ref{f2}), (\ref{sym-E}) and (\ref{sym-E1.17}), we could
obtain the following lemma easily.
\begin{lem}\label{sym-L2.1}
Under the assumptions of Theorem \ref{sym-thm}, we have
    \begin{equation}
      r_t=u,\label{sym-E2.1-00}
    \end{equation}
    \begin{equation}
      A\rho_\infty^\gamma(x)=P_\infty+\int^M_x\frac{Gy}{r_\infty^{2n-2}(y)}dy,
      \label{sym-rhoinf1}
    \end{equation}
    \begin{equation}
  \left(\frac{P_\infty}{A}
  \right)^\frac{1}{\gamma}\leq \rho_\infty \leq
    \bar{\rho}<\infty,
  \ r_\infty^n(x)\in[C^{-1}x,Cx], \label{sym-rhoinf}
  \end{equation}
        \begin{equation}
      \frac{d}{dx}\left(A\rho_\infty^\gamma(x)
      \right)=-G\frac{x}{r_\infty^{2n-2}},\label{sym-rhoinf3}
    \end{equation}
for all $x\in[0,M]$.
\end{lem}

\begin{lem}
Under the assumptions of Theorem \ref{sym-thm}, we have
    \begin{eqnarray}
      &&\frac{d}{dt}\int^M_0\left(\frac{1}{2}u^2+\frac{A\rho^{\gamma-1}}{\gamma-1}
      +\frac{P_\infty}{\rho}+
      \int^r_1G\frac{x}{s^{n-1}}ds\right)dx\nonumber\\
              &&+\int^M_0\{
              (\frac{2}{n}c_1+c_2)\rho^{1+\theta}[(r^{n-1}u)_x]^2+\frac{2(n-1)}{n}
              c_1\rho^{1+\theta}(r^{n-1}u_x-\frac{u}{r\rho})^2\}dx\nonumber\\
      &=&-\int^M_0\Delta f
      udx-\Delta P(ur^{n-1})(M,t).
       \label{sym-E2.1-0}
    \end{eqnarray}
\end{lem}
\begin{proof}
  Multiplying (\ref{sym-E})$_2$ by $u$, integrating the
resulting equation over $[0,M]$, using integration by parts and
 the boundary conditions
(\ref{sym-Efixbd})-(\ref{sym-Efd}), we obtain
    \begin{eqnarray}
      &&\frac{d}{dt}\int^M_0\frac{1}{2}u^2dx-\int^M_0
            A\rho^\gamma\partial_x(r^{n-1}u)dx\nonumber\\
                &&+\int^M_0\left\{
      (2c_1+c_2)\rho^{1+\theta}[(r^{n-1}u)_x]^2-
      2c_1(n-1)\rho^\theta(r^{n-2}u^2)_x\right\}dx\nonumber\\
            &=&
            -P_{\Gamma}(ur^{n-1})(M,t)-\int^M_0     fudx.
            \label{sym-E2.2}
    \end{eqnarray}
From (\ref{sym-E}), we have
    \begin{equation}
   - \int^M_0            A\rho^\gamma\partial_x(r^{n-1}u)dx
            = \frac{d}{dt}\int^M_0\frac{A}{\gamma-1}\rho^{\gamma-1}dx,
    \end{equation}
        \begin{eqnarray}
    &&-P_{\Gamma}(ur^{n-1})(M,t)=
    -P_\infty(r_tr^{n-1})(M,t)
    -\Delta P(ur^{n-1})(M,t)\nonumber\\
            &=&-\frac{d}{dt}\left\{P_\infty\frac{r^n(M,t)}{n}\right\}
            -\Delta P(ur^{n-1})(M,t)\nonumber\\
        &=&-\frac{d}{dt}\int^M_0\frac{P_\infty}{\rho}dx
       -\Delta P(ur^{n-1})(M,t),
        \end{eqnarray}
    \begin{equation}
      -\int^M_0fudx=-\frac{d}{dt}\int^M_0\int^r_1G\frac{x}{s^{n-1}}dsdx-\int^M_0\Delta
      fudx,
    \end{equation}
and
    \begin{eqnarray}
      &&
      (2c_1+c_2)\rho^{1+\theta}(r^{n-1}u)_x^2-
      2c_1(n-1)\rho^\theta(r^{n-2}u^2)_x\nonumber\\
            &=&
      (\frac{2}{n}c_1+c_2)\rho^{1+\theta}(r^{n-1}u)_x^2+\frac{2(n-1)}{n}
      c_1\rho^{1+\theta}(r^{n-1}u_x-\frac{u}{r\rho})^2.\label{sym-E2.5}
    \end{eqnarray}
From (\ref{sym-E2.2})-(\ref{sym-E2.5}), we obtain
 (\ref{sym-E2.1-0}) immediately.
\end{proof}

\textbf{Claim 1}:  Under the assumptions of Theorem \ref{sym-thm},
there is a small positive constant $\epsilon_1$, such that, for
any $T>0$, if
    \begin{equation}
      I(t)=\|\rho(\cdot,t)-\rho_\infty\|_{L^\infty}
      +\left\|\frac{u}{r}(\cdot,t)\right\|_{L^\infty}\leq
      2\epsilon_1,\ \forall\ t\in[0,T],\label{sym-E2.7}
    \end{equation}
then
    $$
    I(t)\leq \epsilon_1, \ \forall\ t\in[0,T].
    $$

    Using the results in Lemmas \ref{sym-L2.2}-\ref{sym-L2.8}, we
    can give the definition of $\epsilon_1$ in
    (\ref{sym-self-E3.76}) and finish the proof of \textbf{Claim
    1}.
\begin{lem}\label{sym-L2.2}
  Under the assumptions of Theorem \ref{sym-thm} and
  (\ref{sym-E2.7}), if $\epsilon_1$ is small enough, we obtain
        \begin{equation}
          \rho(x,t)\in
          \left[\frac{1}{2}\underline{\rho},\frac{3}{2}\bar{\rho}          \right],
          \label{sym-E2.9}
        \end{equation}
  \begin{equation}
    r^n(x,t)\in [C^{-1}x,Cx],
              \label{sym-E2.10}
  \end{equation}
        \begin{equation}
       \|u(\cdot,t)\|_{L^2}+   \|\rho(\cdot,t)-\rho_\infty\|_{L^2}
       +\|r_\infty^{-n}(r^n-r_\infty^n)\|_{L^2_x}\leq C_1\epsilon_0,
                    \label{sym-E2.11}
        \end{equation}
        \begin{equation}
        \int^t_0\int^M_0\left(u^2
        +r^{2n-2}u_x^2+\frac{u^2}{r^2}
        \right)dxds\leq C_1\epsilon_0^2,
        \label{sym-E2.13}
        \end{equation}
for all $t\in[0,T]$ and $x\in[0,M]$.
\end{lem}
\begin{proof}
  From Lemma \ref{sym-L2.1} and (\ref{sym-E2.7}), we can easily obtain the estimate (\ref{sym-E2.9}) when
$2\epsilon_1\leq \frac{1}{2}\underline{\rho}$. From
(\ref{sym-E})$_3$ and (\ref{sym-E2.9}), we can obtain
(\ref{sym-E2.10}) immediately. From (\ref{sym-self-E2.9-1}),
(\ref{sym-rhoinf1}) and (\ref{sym-E2.1-0}), we have
    \begin{eqnarray}
      &&\frac{d}{dt}\left(\int^M_0\frac{1}{2}u^2dx+S[V]-S[V_\infty]\right)\nonumber\\
              &&+\int^M_0\left\{
              (\frac{2}{n}c_1+c_2)\rho^{1+\theta}(r^{n-1}u)_x^2+\frac{2(n-1)}{n}
              c_1\rho^{1+\theta}(r^{n-1}u_x-\frac{u}{r\rho})^2\right\}dx\nonumber\\
      &=&
    -\int^M_0\Delta f udx-\Delta P(ur^{n-1})(M,t)
     \label{sym-E2.16-1}
    \end{eqnarray}
where $V_\infty=\frac{r_\infty^n}{n}$ and $V=\frac{r^n}{n}$. From
 (\ref{epsilon1}),
(\ref{sym-self-E2.9}), (\ref{sym-E2.9})-(\ref{sym-E2.10}) and
Proposition \ref{sym-self-stat-ex-prop3}, we have
        \begin{eqnarray}
      &&C^{-1}\int^M_0\left[(\rho-\rho_\infty)^2+\frac{(V-V_\infty)^2}{V_\infty^{2}}\right]dx\nonumber\\
      &\leq&  S[V]-S[V_\infty]\leq
      C\int^M_0\left[(\rho-\rho_\infty)^2+\frac{(V-V_\infty)^2}{V_\infty^{2}}\right]dx,
        \end{eqnarray}
when $\|V-V_\infty\|_{C^1([0,M])}\leq C_2\epsilon_1\leq \delta_5$,
and
        \begin{equation}
          \left|\Delta P(ur^{n-1})(M,t)\right|\leq C\epsilon_0\left(\int^M_0|\partial_x(r^{n-1}u)|^2dx
          \right)^\frac{1}{2}.\label{sym-E2.18-1}
        \end{equation}
From (\ref{epsilon1})-(\ref{epsilon2}), (\ref{sym-E2.9}) and
(\ref{sym-E2.16-1})-(\ref{sym-E2.18-1}), we obtain
    \begin{eqnarray}
      &&\int^M_0\left(u^2+(\rho-\rho_\infty)^2+r_\infty^{-2n}(r^n-r^n_\infty)^2\right)dx\nonumber\\
              &&+\int^t_0\int^M_0\left\{(r^{n-1}u)_x^2+
              (r^{n-1}u_x-\frac{u}{r\rho})^2\right\}dxds\nonumber\\
      &\leq&C\epsilon_0^2+C\int^t_0 f_1(s)\|u(\cdot,s)\|_{L^2}ds,
    \end{eqnarray}
using Gronwall's inequality
 and (\ref{epsilon1}), we
can obtain (\ref{sym-E2.11})-(\ref{sym-E2.13}) immediately.
\end{proof}

\begin{lem}
  Under the assumptions of Lemma \ref{sym-L2.2}, if $\epsilon_0$ is small enough, we obtain
        \begin{equation}
    \int^t_0\int^M_0\left[(\rho-\rho_\infty)^2+r_\infty^{-2n}(r^n-r_\infty^{n})^2\right]dxds\leq
    C_3\epsilon_0^2,\label{sym-E2.21}
        \end{equation}
for all $t\in[0,T]$.
\end{lem}
\begin{proof}
  Multiplying (\ref{sym-E})$_2$ by $r^{1-n}(\frac{r^n}{n}-\frac{r^n_\infty}{n})$, integrating the
resulting equation over $[0,M]$, using integration by parts and
 the boundary conditions
(\ref{sym-Efixbd})-(\ref{sym-Efd}), we obtain
    \begin{eqnarray}
      &&\int^M_0A(\rho_\infty^\gamma-\rho^\gamma)(\rho^{-1}-\rho_\infty^{-1})
      +Gx(r^{2-2n}-r_\infty^{2-2n})(\frac{r^n}{n}-\frac{r^n_\infty}{n})dx\nonumber\\
            &=&-\int^M_0\frac{u_t}{r^{n-1}}(\frac{r^n}{n}-\frac{r^n_\infty}{n})dx
            +\Delta P\left.\left\{\frac{r^n}{n}-\frac{r^n_\infty}{n}\right\}\right|_{x=M}
            \nonumber\\
            &&-\int^M_0\Delta f\frac{r^{1-n}}{n}(r^n-r^n_\infty
      )dx+\int^M_0
            2c_1(n-1)\rho^\theta \left(\frac{u}{r}(\frac{r^n}{n}-\frac{r^n_\infty}{n})
                   \right)_xdx\nonumber\\
     &&+\int^M_0
            (2c_1+c_2)\rho^{1+\theta}\partial_x(r^{n-1}u)(\rho_\infty^{-1}-\rho^{-1})dx\nonumber\\
    &:=&\sum^5_{i=1}I_i.\label{sym-E2.22-1}
    \end{eqnarray}
We can rewrite the left hand side of (\ref{sym-E2.22-1}) as
follows
    \begin{eqnarray*}
      &&\textrm{L.H.S of } (\ref{sym-E2.22-1})\nonumber\\&=&\int^M_0
      \left[\gamma
      A\rho_\infty^{\gamma+1}(\rho^{-1}-\rho_\infty^{-1})^2-(2n-2)Gxr_\infty^{2-3n}
      \left(\frac{r^n}{n}-\frac{r_\infty^n}{n}\right)^2
      \right]dx\nonumber\\
            &&+\int^M_0
      \left[g_1(\rho^{-1}-\rho_\infty^{-1})^2+g_2r_\infty^{-2n}
      \left(\frac{r^n}{n}-\frac{r_\infty^n}{n}\right)^2
      \right]dx,
    \end{eqnarray*}
where
    $$
    |g_1|=\left|\frac{A(\rho_\infty^\gamma-\rho^\gamma)}{\rho^{-1}-\rho_\infty^{-1}}-\gamma
    A\rho_\infty^{1+\gamma}\right|\leq C_4\epsilon_1
    $$
and
    $$
    |g_2|=\left|Gxr_\infty^{2n}(r^{2-2n}-r_\infty^{2-2n})\left(\frac{r^n}{n}-\frac{r_\infty^n}{n}\right)^{-1}
    +(2n-2)Gxr_\infty^{2-n}
    \right|\leq C_4\epsilon_1.
    $$
From (\ref{sym-selfE2.8}), we have
    \begin{eqnarray}
    \textrm{L.H.S of } (\ref{sym-E2.22-1})&\geq&
    (2C_5-C_4\epsilon_1)\int^M_0
      \left[(\rho^{-1}-\rho_\infty^{-1})^2+r_\infty^{-2n}
      (\frac{r^n}{n}-\frac{r_\infty^n}{n})^2
      \right]dx\nonumber\\
                &\geq&C_5\int^M_0
      \left[(\rho^{-1}-\rho_\infty^{-1})^2+r_\infty^{-2n}
      \left(\frac{r^n}{n}-\frac{r_\infty^n}{n}\right)^2
      \right]dx,
    \end{eqnarray}
when $C_4\epsilon_1\leq C_5$.

 From (\ref{sym-E2.7}) and (\ref{sym-E2.9})-(\ref{sym-E2.10}),
using integration by parts, we can estimate $I_i$ as follows.
\begin{eqnarray}
  I_1&=&-\frac{d}{dt}\int^M_0\frac{u}{r^{n-1}}\left(\frac{r^n}{n}-\frac{r_\infty^n}{n}
      \right)dx+\int^M_0u^2\left(\frac{1}{n}+\frac{(n-1)r_\infty^n}{nr^n}
      \right)dx\nonumber\\
            &\leq& -\frac{d}{dt}\int^M_0\frac{u}{r^{n-1}}\left(\frac{r^n}{n}-\frac{r_\infty^n}{n}
      \right)dx+ C\int^M_0u^2dx,
\end{eqnarray}
        \begin{equation}
          I_2=\Delta P\int^M_0(\rho^{-1}-\rho_\infty^{-1})dx\leq
          \frac{C_5}{10}\int^M_0(\rho^{-1}-\rho_\infty^{-1})^2dx+C|\Delta P|^2,
        \end{equation}
    \begin{equation}
      I_3\leq\frac{C_5}{10}\int^M_0r_\infty^{-2n}
      \left(\frac{r^n}{n}-\frac{r_\infty^n}{n}\right)^2dx+Cf_1^2,
    \end{equation}
    \begin{eqnarray}
      I_4&\leq&\frac{C_5}{10}\int^M_0\left[(\rho^{-1}-\rho_\infty^{-1})^2+r_\infty^{-2n}
       \left(\frac{r^n}{n}-\frac{r_\infty^n}{n}\right)^2
      \right]dx\nonumber\\
      &&+C\int^M_0\left(
      [(r^{n-1}u)_x]^2+\frac{u^2}{r^2}
      \right)dx
    \end{eqnarray}
and
    \begin{equation}
      I_5\leq\frac{C_5}{10}\int^M_0(\rho^{-1}-\rho_\infty^{-1})^2dx+C\int^M_0
      (r^{n-1}u)_x^2dx.\label{sym-E2.30-1}
    \end{equation}
From (\ref{sym-E2.22-1})-(\ref{sym-E2.30-1}), we get
    \begin{eqnarray}
    &&\frac{d}{dt}\int^M_0\frac{u}{nr^{n-1}}(r^n-r_\infty^n
      )dx+C\int^M_0[(\rho^{-1}-\rho_\infty^{-1})^2+r_\infty^{-2n}
       (\frac{r^n}{n}-\frac{r_\infty^n}{n})^2
      ]dx\nonumber\\
      &  \leq& C\int^M_0\left(r^{2n-2}u_x^2+\frac{u^2}{r^2}\right)dx
      +C\left(|\Delta P|^2+f_1^2
          \right).\label{sym-E2.33-1}
    \end{eqnarray}
And from (\ref{sym-E2.9})-(\ref{sym-E2.13}), we obtain
(\ref{sym-E2.21}) immediately.
\end{proof}

From now on, we study the case of $\theta>0$ and prove the main
results in this case only, since the case of $\theta=0$ can be
discussed through the similar process.

\begin{lem}\label{sym-L2.3}
  Under the assumptions of  Lemma \ref{sym-L2.2},
  if $\epsilon_1$  is small enough, we  obtain
        \begin{equation}
          \int^M_0[r^{2n-2}(\rho-\rho_\infty)_x^2](x,t)dx+\int^t_0\int^M_0
          [r^{2n-2}(\rho-\rho_\infty)_x^2](x,s)dxds\leq C_6
          \epsilon_0^2,
          \label{sym-E2.14}
        \end{equation}
    for all $t\in[0,T]$.
\end{lem}
\begin{proof}
  From (\ref{sym-E}), we have
    \begin{eqnarray}
      &&\partial_tH+\frac{A\gamma\rho^{\gamma-\theta}}{2c_1+c_2}H\nonumber\\
            &=&\frac{A\gamma}{2c_1+c_2}\rho^{\gamma-\theta}u+
            (\frac{2c_1+c_2}{\theta}-2c_1)(n-1)r^{n-2}u(\rho^\theta)_x-f(x,r,t)\nonumber\\
                &&-\frac{(n-1)(2c_1+c_2)}{\theta}r^{n-2}u(\rho_\infty^\theta)_x-\frac{A\gamma \rho^{\gamma-\theta}r^{n-1}}{\theta}
                (\rho_\infty^\theta)_x.
            \label{sym-E2.15}
    \end{eqnarray}
where
$H=u+\frac{2c_1+c_2}{\theta}r^{n-1}(\rho^\theta-\rho_\infty^\theta)_x
$. Multiplying (\ref{sym-E2.15}) by $H$, integrating the resulting
equation over $[0,M]$, using the Cauchy-Schwarz inequality, we
obtain
    \begin{eqnarray}
      &&\frac{d}{dt}\int^M_0H^2(x,t)dx+C_7\int^M_0H^2(x,t)dx\nonumber\\
      &\leq&C\int^M_0\left(
      |H\rho^{\gamma-\theta}u|+|\frac{u}{r}H^2|
      +|\frac{u^2}{r}H|+|\Delta fH|\right)dx\nonumber\\
                &&+\int^M_0\left|G\frac{x}{r^{n-1}}+\frac{ \rho^{\gamma-\theta}r^{n-1}}{ \rho_\infty^{\gamma-\theta}}
                (A\rho_\infty^\gamma)_x\right||H|
                dx+C\int^M_0|r^{n-2}u(\rho_\infty^\gamma)_xH|
                dx\nonumber\\
            &\leq& C\int^M_0\left(
            u^2+\frac{u^4}{r^2}+\frac{x^2r^{2n-4}}{r_\infty^{4n-4}}u^2
            \right)dx +(\frac{1}{4}+C_8\epsilon_1) C_7\int^M_0H^2dx\nonumber\\
              &&      +C\int^M_0\left|G\frac{x}{r^{n-1}}+\frac{ \rho^{\gamma-\theta}r^{n-1}}{ \rho_\infty^{\gamma-\theta}}
                (A\rho_\infty^\gamma)_x\right|^2dx+Cf_1^2,
            \label{sym-E2.16}
    \end{eqnarray}
From (\ref{sym-rhoinf3}) and (\ref{sym-E2.9}),
 we have
    \begin{eqnarray}
   && \int^M_0\left|G\frac{x}{r^{n-1}}+\frac{ \rho^{\gamma-\theta}r^{n-1}}{ \rho_\infty^{\gamma-\theta}}
                (A\rho_\infty^\gamma)_x\right|^2dx\nonumber\\
                    &=&C\int^M_0\left|G\frac{x}{r^{n-1}}
                    -G\frac{ x\rho^{\gamma-\theta}r^{n-1}}{ \rho_\infty^{\gamma-\theta}r_\infty^{2n-2}}
                \right|^2dx\nonumber\\
        &\leq&
        C\int^M_0\left[(r-r_\infty)^2+(\rho-\rho_\infty)^2
           \right]dx.\label{sym-E2.24-1}
    \end{eqnarray}
Then, if $\epsilon_1\leq1$ and $C_8\epsilon_1\leq \frac{1}{4}$,
from (\ref{sym-E2.16})-(\ref{sym-E2.24-1}), we obtain
    \begin{eqnarray}
      &&\frac{d}{dt}\int^M_0H^2(x,t)dx+\frac{C_7}{2}\int^M_0H^2(x,t)dx\nonumber\\
      &\leq& C\int^M_0\left(
            u^2+(r-r_\infty)^2+(\rho-\rho_\infty)^2
            \right)dx+Cf_1^2.
            \label{sym-E2.36-4}
    \end{eqnarray}
From  (\ref{sym-E2.9})-(\ref{sym-E2.13}), (\ref{sym-E2.21}) and
(\ref{sym-E2.36-4}), we obtain (\ref{sym-E2.14}) immediately.
\end{proof}

\begin{lem}
   Under the assumptions of Lemma \ref{sym-L2.2}, if $\epsilon_1$ is small enough, we obtain
        \begin{equation}
         (\rho(M,t)-\rho_\infty(M))^2+ \int^t_0(\rho(M,s)-\rho_\infty(M))^2ds\leq C_{10}\epsilon_0^2,
         \label{sym-E2.24}
        \end{equation}
        \begin{equation}
          \int^t_0\int^M_0r^{\frac{1}{2}-m}(\rho-\rho_\infty)^2dxds\leq
          C_{11}\epsilon_0^2,\label{sym-E2.25}
        \end{equation}
        \begin{equation}
         \int^M_0(r^{\frac{1}{2}-m}
         u^2)(x,t)dx+\int^t_0\int^M_0r^{\frac{1}{2}-m}\left(
         r^{2n-2}u_x^2+\frac{u^2}{r^2}
         \right)dxds\leq C_{12}\epsilon_0^2,\label{sym-E2.27}
        \end{equation}
        \begin{equation}
          \int^M_0(r^{2n-2+{\frac{1}{2}-m}}(\rho-\rho_\infty)_x^2)dx+\int^t_0\int^M_0
          (r^{2n-2+{\frac{1}{2}-m}}(\rho-\rho_\infty)_x^2)dxds\leq C_{15}
          \epsilon_0^2,
          \label{sym-E2.32}
        \end{equation}
                \begin{equation}
          \|\rho(\cdot,t)-\rho_\infty(\cdot)\|_{L^\infty}
          +\int^M_0|(\rho-\rho_\infty)_x|(x,t)dx\leq C_{16}
          \epsilon_0,
          \label{sym-E2.32-1}
        \end{equation}
            \begin{equation}
              |r(x,t)-r_\infty(x)|\leq
              C_{16}\epsilon_0x^{\frac{1}{n}},   \ x\in[0,M],       \label{sym-E2.32-1-1}
            \end{equation}
for all $t\in[0,T]$ and $m=0,1,\ldots,n-1$.
\end{lem}
\begin{proof}
  From (\ref{sym-E})$_1$ and the boundary condition (\ref{sym-Efd}), we have
    $$
    A\rho^\gamma(M,t)-P_\Gamma+\frac{(2c_1+c_2)}{\theta}\partial_t(\rho^\theta)(M,t)
    =-2c_1(n-1)\left(\rho^\theta\frac{u}{r}\right)(M,t).
    $$
  Multiplying the above equality by $\rho^\theta(M,t)-\rho_\infty^\theta(M)$,   we obtain
        \begin{eqnarray*}
        &&\left.\frac{2c_1+c_2}{2\theta}\frac{d}{dt}(\rho^\theta-\rho_\infty^\theta)^2\right|_{x=M}
    +(\rho^\theta(M,t)-\rho_\infty^\theta(M))(A\rho^\gamma(M,t)-P_\infty)\nonumber\\
            &=& -2c_1(n-1)\left.\left[\rho^\theta\frac{u}{r}(\rho^\theta-\rho_\infty^\theta)
            \right]\right|_{x=M}   +\Delta P(\rho^\theta(M,t)-\rho_\infty^\theta(M)).
        \end{eqnarray*}
Combining (\ref{sym-E2.9})-(\ref{sym-E2.10}), using the
Cauchy-Schwarz inequality, we get
    \begin{eqnarray}
      &&\frac{d}{dt}(\rho^\theta(M,t)-\rho^\theta_\infty(M))^2+C^{-1}(\rho^\theta(M,t)-\rho^\theta_\infty(M))^2\nonumber\\
            &\leq&C|\Delta P|^2+C(u^2r^n)(M,t)=C|\Delta P|^2+C\int^M_0\partial_x(u^2r^n)dx\nonumber\\
            &\leq&C|\Delta
            P|^2+C\int^M_0(r^{2n-2}u_x^2+\frac{u^2}{r^2})dx.\label{sym-E2.53-1}
    \end{eqnarray}
Integrating the above inequality over $[0,t]$, using the estimates
(\ref{sym-E2.9}) and (\ref{sym-E2.13}), we can obtain
(\ref{sym-E2.24}).

From (\ref{sym-E2.11})-(\ref{sym-E2.13}), (\ref{sym-E2.21}) and
(\ref{sym-E2.14}), we know that the estimates
(\ref{sym-E2.25})-(\ref{sym-E2.32}) hold with $m=0$.

\textbf{Claim 2}: If that (\ref{sym-E2.25})-(\ref{sym-E2.32}) hold
with $m\leq k$, $k\in [0,n-2]$, then the estimates
(\ref{sym-E2.25})-(\ref{sym-E2.32}) hold with $m= k+1$.

We could prove Claim 2 as follows. Let $\alpha_k=\frac{1}{2}-k-1$.
Using H\"{o}lder's inequality, we have
    \begin{eqnarray}
      &&\int^M_0r^{\alpha_k}(\rho-\rho_\infty)^2dx\nonumber\\
      &=&
      \int^M_0r^{\alpha_k}\left(\rho(M,s)-\rho_\infty(M)-\int^M_x\partial_x
     ( \rho-\rho_\infty)dy
      \right)^2dx\nonumber\\
      &\leq&C(\rho(M,s)-\rho_\infty(M))^2\nonumber\\
            &&+C\int^M_0r^{\alpha_k}\int^M_x r^{2n-2+\alpha_{k-1}}
      (\rho-\rho_\infty)^2_xdy
      \int^M_xr^{2-2n-\alpha_{k-1}}dy
      dx\nonumber\\
                &\leq&C(\rho(M,s)-\rho_\infty(M))^2+C\int^M_0 r^{2n-2+\alpha_{k-1}}
      (\rho-\rho_\infty)^2_xdx.\label{sym-E2.53}
    \end{eqnarray}
From (\ref{sym-E2.24}), (\ref{sym-E2.32})($m=k$) and
(\ref{sym-E2.53}), we obtain (\ref{sym-E2.25}) ($m=k+1$).

  Multiplying (\ref{sym-E})$_2$ by $ur^{\alpha_k}$, integrating the
resulting equation over $[0,M]$, using integration by parts and
the boundary conditions (\ref{sym-Efixbd})-(\ref{sym-Efd}), we
obtain
    \begin{eqnarray}
      &&\frac{d}{dt}\int^M_0\frac{1}{2}r^{\alpha_k} u^2dx
      -\int^M_0\frac{{\alpha_k}}{2}r^{{\alpha_k}-1}u^3dx\nonumber\\
                &=&-\int^M_0[
                (2c_1+c_2)\rho^{1+\theta} (r^{n-1}u)_x(r^{n-1+{\alpha_k}}u)_x
                -2c_1(n-1)\rho^\theta(r^{n-2+{\alpha_k}}u^2)_x
                ]dx\nonumber\\
      &&+\int^M_0A(\rho^\gamma-\rho_\infty^\gamma)(r^{n-1+{\alpha_k}}u)_xdx
      -\int^M_0\Delta fur^{\alpha_k} dx\nonumber\\
            &&  -\Delta P
        (ur^{n-1+{\alpha_k}})(M,t)      +\int^M_0(r^{n-1+{\alpha_k}}u)_x\int^M_x\left(
        \frac{Gy}{r_\infty^{2n-2}}-\frac{Gy}{r^{2n-2}}
        \right)dydx
      \nonumber\\
        &&\nonumber\\
      &=&\sum^5_{i=1}L_i.\label{sym-E2.28}
    \end{eqnarray}
We can estimate $L_1$  as follows.
    \begin{eqnarray}
      -L_1&=&\int^M_0\left\{(2c_1+c_2)r^{2n-2+{\alpha_k}}\rho^{1+\theta}u_x^2
      \right.\nonumber\\
            &&+[2{\alpha_k} c_1+c_2(2n-2+{\alpha_k})]\rho^\theta r^{n-2+{\alpha_k}}uu_x
        \nonumber\\
      &&\left.
      + [2c_1(n-1)+c_2(n-1)(n-1+{\alpha_k})]\rho^{\theta-1}r^{{\alpha_k}-2}u^2
                \right\}dx.
    \end{eqnarray}
Since
$$[2{\alpha}
c_1+c_2(2n-2+{\alpha})]^2-4(2c_1+c_2)[2c_1(n-1)+c_2(n-1)(n-1+{\alpha})]<0,$$
where $\alpha=\frac{3}{2}-n$ and
$$[c_2(2n-2)]^2-4(2c_1+c_2)[2c_1(n-1)+c_2(n-1)^2]<0,$$
we have
$$[2{\alpha_k}
c_1+c_2(2n-2+{\alpha_k})]^2-4(2c_1+c_2)[2c_1(n-1)+c_2(n-1)(n-1+{\alpha_k})]<0.$$
Then there exists a positive constant $C_{13}$ such that
    \begin{equation}
      -L_1\geq C_{13}\int^M_0\left(r^{2n-2+{\alpha_k}}\rho^{1+\theta}u_x^2+\rho^{\theta-1}r^{{\alpha_k}-2}u^2
      \right)dx.\label{sym-E2.30}
    \end{equation}
From (\ref{sym-E2.9})-(\ref{sym-E2.10}), using the Cauchy-Schwarz
inequality, we obtain
    \begin{eqnarray}
      L_2&=&\int^M_0A(\rho^\gamma-\rho_\infty^\gamma)r^{n-1+{\alpha_k}}u_xdx
      +\int^M_0A(n-1+{\alpha_k})r^{{\alpha_k}-1}u\frac{\rho^\gamma-\rho_\infty^\gamma}{\rho}dx
      \nonumber\\
            &\leq& \frac{C_{13}}{8}\int^M_0\left(r^{2n-2+{\alpha_k}}\rho^{1+\theta}u_x^2+\rho^{\theta-1}
            r^{{\alpha_k}-2}u^2
            \right)dx+C\int^M_0r^{\alpha_k}(\rho-\rho_\infty)^2dx,
      \label{sym-E2.31}
    \end{eqnarray}
    \begin{equation}
      L_3  \leq\frac{1}{8}C_{13}\int^M_0\rho^{\theta-1}r^{{\alpha_k}-2}u^2
      dx+Cf_1^2,
    \end{equation}
        \begin{eqnarray}
      L_4  &=&-\int^M_0\Delta P\partial_x(r^{n-1+{\alpha_k}}u)dx\nonumber\\
                &\leq& \frac{1}{8}C_{13}\int^M_0\left(r^{2n-2+{\alpha_k}}
                \rho^{1+\theta}u_x^2+\rho^{\theta-1}r^{{\alpha_k}-2}u^2
                \right)dx+C|\Delta P|^2
    \end{eqnarray}
and
    \begin{eqnarray}
      L_5&\leq& \frac{1}{8}C_{13}\int^M_0\left(r^{2n-2+{\alpha_k}}\rho^{1+\theta}u_x^2+
      \rho^{\theta-1}r^{{\alpha_k}-2}u^2      \right)dx\nonumber\\
            &&+C\int^M_0r^{\alpha_k}\left[\int^M_x\left(
        \frac{Gy}{r_\infty^{2n-2}}-\frac{Gy}{r^{2n-2}}
        \right)dy\right]^2dx\nonumber\\
                           &\leq& \frac{C_{13}}{8}\int^M_0\left(r^{2n-2+{\alpha_k}}
                           \rho^{1+\theta}u_x^2+\rho^{\theta-1}r^{{\alpha_k}-2}u^2
                            \right)dx+C\int^M_0r^{\alpha_k}(\rho-\rho_\infty)^2dx.\label{sym-E2.56}
    \end{eqnarray}
 From (\ref{sym-E2.7}), (\ref{sym-E2.9})-(\ref{sym-E2.10})
 and
(\ref{sym-E2.28})-(\ref{sym-E2.56}), we obtain
        \begin{eqnarray*}
         &&\frac{d}{dt}\int^M_0(r^{\alpha_k}
         u^2)(x,t)dx+2C_{14}^{-1}\int^M_0r^{\alpha_k}\left(
         r^{2n-2}u_x^2+\frac{u^2}{r^2}
         \right)dx\\
                &\leq& C(f_1^2+|\Delta P|^2)+C_{14}\epsilon_1\int^M_0r^{\alpha_k}
         \frac{u^2}{r^2} dx+C\int^M_0r^{\alpha_k}(\rho-\rho_\infty)^2dx.
        \end{eqnarray*}
When $C_{14}^2\epsilon_1\leq \frac{1}{2}$, using the estimate
(\ref{sym-E2.25}) ($m=k+1$), we can get
    \begin{eqnarray}
    &&\frac{d}{dt}\int^M_0(r^{\alpha_k}
         u^2)(x,t)dx+C_{14}^{-1}\int^M_0r^{\alpha_k}\left(
         r^{2n-2}u_x^2+\frac{u^2}{r^2}
         \right)dx\nonumber\\
         &\leq& C(f_1^2+|\Delta P|^2)+C\int^M_0r^{\alpha_k}(\rho-\rho_\infty)^2dx\label{sym-E2.63}
    \end{eqnarray}
and (\ref{sym-E2.27}) ($m=k+1$) holds.

  From (\ref{sym-E}), we have
    \begin{eqnarray}
      &&\partial_tH_1+\frac{A\gamma}{2c_1+c_2}\rho^{\gamma-\theta}H_1\nonumber\\
            &=&\frac{{\alpha_k}}{2}r^{\frac{{\alpha_k}}{2}-1}u^2+\frac{A\gamma}{2c_1+c_2}\rho^{\gamma-\theta}r^\frac{{\alpha_k}}{2}u
            -\frac{A\gamma \rho^{\gamma-\theta}r^{n-1+\frac{{\alpha_k}}{2}}}{\theta}
                (\rho_\infty^\theta)_x\nonumber\\
                &&+ \left[\frac{2c_1+c_2}{\theta}(n-1+\frac{{\alpha_k}}{2})-2c_1(n-1)
            \right]r^{n-2+\frac{{\alpha_k}}{2}}u(\rho^\theta)_x
            \nonumber\\
                    &&-\frac{(n-1+\frac{{\alpha_k}}{2})(2c_1+c_2)}{\theta}r^{n-2+\frac{{\alpha_k}}{2}
                    }u(\rho_\infty^\theta)_x-f(x,r,t)r^{\frac{{\alpha_k}}{2}}.
            \label{sym-E2.33}
    \end{eqnarray}
where
$H_1=r^\frac{{\alpha_k}}{2}u+\frac{2c_1+c_2}{\theta}r^{n-1+\frac{{\alpha_k}}{2}}(\rho^\theta-\rho_\infty^\theta)_x$.
Multiplying (\ref{sym-E2.33}) by $H_1$, integrating the resulting
equation over $[0,M]$, using the Cauchy-Schwarz inequality, we
obtain
    \begin{eqnarray}
      &&\frac{d}{dt}\int^M_0H_1^2(x,t)dx+C_{18}\int^M_0H_1^2(x,t)dx\nonumber\\
      &\leq&
      C\int^M_0\left(
      |H_1\rho^{\gamma-\theta}r^\frac{{\alpha_k}}{2}u|
      +\left|\frac{u}{r}\right|H_1^2
      +r^{\frac{{\alpha_k}}{2}-1}u^2|H_1|+r^\frac{{\alpha_k}}{2}|\Delta fH_1|\right)dx\nonumber\\
            &&+C\int^M_0\left|r^{n-2+\frac{{\alpha_k}}{2}}
            u(\rho_\infty^\gamma)_xH_1
            \right|+r^{\frac{{\alpha_k}}{2}}
            \left|\frac{Gx}{r^{n-1}}+\frac{\rho^{\gamma-\theta}r^{n-1}}{\rho_\infty^{\gamma-\theta}}(A\rho_\infty^\gamma)_x
            \right||H_1|dx\nonumber\\
            &\leq& C\int^M_0\left(
            r^{\alpha_k}
            u^2+r^{{\alpha_k}-2}u^4+\frac{x^2r^{2n-4+{\alpha_k}}}{r_\infty^{4n-4}}u^2
            \right)dx+Cf_1^2\nonumber\\
                    &&+(\frac{1}{4}+C_{17}\epsilon_1)C_{18}\int^M_0H_1^2dx+C\int^M_0r^{{\alpha_k}}
            \left|\frac{Gx}{r^{n-1}}+\frac{\rho^{\gamma-\theta}r^{n-1}}{\rho_\infty^{\gamma-\theta}}(\rho_\infty^\gamma)_x
            \right|^2dx.
            \label{sym-E2.34}
    \end{eqnarray}
If $\epsilon_1\leq 1,C_{17}\epsilon_1\leq \frac{1}{4}$, using the
estimates (\ref{sym-E2.7}) and (\ref{sym-E2.9})-(\ref{sym-E2.10}),
we have
    \begin{eqnarray}
    &&\frac{d}{dt}\int^M_0H_1^2(x,t)dx+\frac{C_{18}}{2}\int^M_0H_1^2(x,t)dx\nonumber\\
    &\leq&
    C\int^M_0\left[r^{\alpha_k}u^2+r^{\alpha_k}(\rho-\rho_\infty)^2
    +r^{\alpha_k}(r-r_\infty)^2
    \right]dx+Cf_1^2\nonumber\\
    &\leq&
    C\int^M_0\left[r^{\alpha_k}u^2+r^{\alpha_k}(\rho-\rho_\infty)^2
    \right]dx+Cf_1^2.\label{sym-E2.66}
    \end{eqnarray}
And from (\ref{sym-E2.25})-(\ref{sym-E2.27}) ($m=k+1$), we have
    $$
    \int^M_0H_1^2(x,t)dx+\int^t_0\int^M_0H_1^2(x,s)dxds\leq
    C\epsilon_0^2.
    $$

From (\ref{sym-E2.9}) and (\ref{sym-E2.25})-(\ref{sym-E2.27})
($m=k+1$), we obtain (\ref{sym-E2.32}) ($m=k+1$) immediately  and
finish the proof of \textbf{Claim 2}.

From Claim 2, we obtain that the estimates
(\ref{sym-E2.25})-(\ref{sym-E2.32}) ($m=0,\ldots,n-1$) hold.
 From (\ref{sym-E2.10}) and
(\ref{sym-E2.32}), using H\"{o}lder's inequality, we obtain
    \begin{equation}
    \int^M_0|(\rho-\rho_\infty)_x|dx
    \leq \left(\int^M_0r^{2n-2+{\alpha}}
    (\rho-\rho_\infty)_x^2dx
    \right)^\frac{1}{2}\left(\int^M_0r^{-2n+2-{\alpha}}dx
    \right)^\frac{1}{2}
    \leq C\epsilon_0.\label{sym-E2.61}
    \end{equation}
From (\ref{sym-E2.10})-(\ref{sym-E2.11}) and (\ref{sym-E2.61}),
using Sobolev's embedding Theorem, we could obtain
(\ref{sym-E2.32-1})-(\ref{sym-E2.32-1-1}) immediately.
\end{proof}

\begin{lem}
 Under the assumptions of  Lemma \ref{sym-L2.2},
  if $\epsilon_1$ is small enough, we obtain
            \begin{equation}
              \int^M_0\left(
              \frac{u^2}{r^2}+r^{2n-2}u_x^2
              \right)(x,t)dx+\int^t_0\int^M_0u_t^2(x,s)dxds\leq
              C_9\epsilon_0^2(1+\|(r^{n-1}u)_x\|_{L^\infty_{tx}}),\label{sym-E2.17}
            \end{equation}
  for all $t\in[0,T]$.
\end{lem}
\begin{proof}
   Multiplying (\ref{sym-E})$_2$ by $u_t$, integrating the
resulting equation over $[0,M]$, using integration by parts and
the boundary conditions (\ref{sym-Efixbd})-(\ref{sym-Efd}), we
obtain
    \begin{eqnarray}
      &&\int^M_0u_t^2 dx+\int^M_0
      (2c_1+c_2)\rho^{1+\theta}(r^{n-1}u)_x
      (r^{n-1}u_t)_xdx\nonumber\\
            &=&\int^M_0A\rho^\gamma(r^{n-1}u_t)_xdx-
            P_{\Gamma}(r^{n-1}u_t)(M,t)\nonumber\\
                &&+\int^M_02c_1(n-1)\rho^\theta(r^{n-2}uu_t)_x
      dx-\int^M_0fu_tdx\nonumber\\
            &:=&\sum^4_{i=1}N_i.\label{sym-E2.18}
    \end{eqnarray}
From (\ref{sym-E2.7}), (\ref{sym-E2.9}) and (\ref{sym-E2.13}),
using the Cauchy-Schwarz inequality, we obtain
    \begin{eqnarray}
      &&\int^M_0
      (2c_1+c_2)\rho^{1+\theta} (r^{n-1}u)_x
       (r^{n-1}u_t)_xdx\nonumber\\
            &=&\frac{d}{dt}\int^M_0\frac{2c_1+c_2}{2}\rho^{1+\theta}[(r^{n-1}u)_x]^2dx\nonumber\\
            &&-\int^M_0
                    (2c_1+c_2)(n-1)\rho^{1+\theta} (r^{n-1}u)_x
                    (r^{n-2}u^2)_xdx\nonumber\\
        &&+\int^M_0
      \frac{(2c_1+c_2)}{2}(1+\theta)\rho^{2+\theta}[ (r^{n-1}u)_x]^3dx\nonumber\\
                        &\geq&\frac{d}{dt}\int^M_0\frac{2c_1+c_2}{2}\rho^{1+\theta}[(r^{n-1}u)_x]^2dx
                        \nonumber\\
                    &&-C(
            \| (r^{n-1}u)_x\|_{L^2_{x}}^2+\left\|\frac{u}{r}\right\|_{L^2_{x}}^2)(1+
            \| (r^{n-1}u)_x\|_{L^\infty_{t,x}}),
    \end{eqnarray}
        \begin{eqnarray}
          N_1 &=&\frac{d}{dt}\int^M_0A\rho^\gamma (r^{n-1}u)_xdx
          +\int^M_0A\gamma\rho^{\gamma+1}[(r^{n-1}u)_x]^2dx\nonumber\\
                        &&-\int^M_02A(n-1)\rho^{\gamma}\frac{u}{r} (r^{n-1}u)_xdx
                        +\int^M_0An(n-1)\rho^{\gamma-1}\frac{u^2}{r^2}dx\nonumber\\
          &\leq&\frac{d}{dt}\int^M_0A\rho^\gamma (r^{n-1}u)_xdx
          +C(\| (r^{n-1}u)_x\|_{L^2_{x}}^2
          +\left\|\frac{u}{r}\right\|_{L^2_{x}}^2),
        \end{eqnarray}
                \begin{eqnarray}
             N_2 &=&-\frac{d}{dt}\int^M_0(P_\infty+\Delta P(t)) (r^{n-1}u)_xdx
               \nonumber\\
                       &&+\int^M_0(n-1)
            (P_\infty+\Delta P)(r^{n-2}u^2)_xdxds
            +(\Delta P)'\int^M_0(r^{n-1}u)_xdyds\nonumber\\
                            &\leq&-\frac{d}{dt}\int^M_0(P_\infty+\Delta P(t)) (r^{n-1}u)_xdx
                            \nonumber\\
            &&                            +C\left(\|(r^{n-1}u)_x\|_{L^2_x}^2+
            \left\|\frac{u}{r}\right\|_{L^2_x}^2+
            |(\Delta   P)'|^2
                            \right),
                 \label{sym-E2.22}
                \end{eqnarray}
    \begin{eqnarray}
   N_3    &=&\frac{d}{dt}\int^M_0c_1(n-1)\rho^\theta(r^{n-2}u^2)_x
                dx+\int^M_02\theta
                c_1(n-1)\rho^{\theta+1}\frac{u}{r}[(r^{n-1}u)_x]^2
                dx\nonumber\\
    &&-\int^M_0\theta
    c_1n(n-1)\rho^\theta\frac{u^2}{r^2} (r^{n-1}u)_xdx+\int^M_0
            2nc_1(n-1)(n-2)\rho^{\theta-1}\frac{u^3}{r^3}dx\nonumber\\
    &&-\int^M_03c_1(n-1)(n-2)\rho^\theta
            \frac{u^2}{r^2} (r^{n-1}u)_xdx\nonumber\\
    &\leq&\frac{d}{dt}\int^M_0c_1(n-1)\rho^\theta(r^{n-2}u^2)_x
              dx+C
    (\| (r^{n-1}u)_x\|_{L^2_{x}}^2+\left\|\frac{u}{r}\right\|_{L^2_{x}}^2)
    \end{eqnarray}
and
        \begin{equation}
          N_4\leq-\frac{d}{dt}\int^M_0G\frac{xu}{r^{n-1}}dx
          +\int^M_0(1-n)Gxr^{-n}u^2dx+\frac{1}{2}\int^M_0u_t^2dx+Cf_1^2.\label{sym-E2.22-2}
        \end{equation}
From (\ref{sym-E2.18})-(\ref{sym-E2.22-2}), using the fact that
    \begin{eqnarray*}
      &&\int^M_0\left\{
      \frac{1}{2}(2c_1+c_2)\rho^{1+\theta}[(r^{n-1}u)_x]^2-
      c_1(n-1)\rho^\theta(r^{n-2}u^2)_x\right\}dx\nonumber\\
            &=&\int^M_0\left\{
      \frac{1}{2}(\frac{2}{n}c_1+c_2)\rho^{1+\theta}[(r^{n-1}u)_x]^2+\frac{(n-1)}{n}
      c_1\rho^\theta(r^{n-1}u_x-\frac{u}{r\rho})^2\right\}dx,
    \end{eqnarray*}
we have
    \begin{eqnarray}
            &&\int^M_0\frac{1}{2}u_t^2dx\nonumber\\
      && +\frac{d}{dt}\int^M_0\left\{
      \frac{1}{2}(\frac{2}{n}c_1+c_2)\rho^{1+\theta}(r^{n-1}u)_x^2+\frac{(n-1)}{n}
      c_1\rho^\theta(r^{n-1}u_x-\frac{u}{r\rho})^2\right\}dx\nonumber\\
            &\leq&\frac{d}{dt}\left\{\int^M_0\left[(A\rho^\gamma-A\rho_\infty^\gamma-\Delta P)
            (r^{n-1}u)_x
            +(r^{n-1}u)_x\int^M_x\left(\frac{Gy}{r_\infty^{2n-2}}-
            \frac{Gy}{r^{2n-2}}
            \right)dy\right]dx\right\}\nonumber\\
                    &&+C(1+\|(r^{n-1}u)_x\|_{L^\infty_{x}})
    (\|(r^{n-1}u)_x\|_{L^2_{x}}^2+\left\|\frac{u}{r}\right\|_{L^2_{x}}^2+f_1^2+|(\Delta P)'|^2).
    \label{sym-E2.45}
    \end{eqnarray}
Integrating (\ref{sym-E2.45}) over $[0,t]$, using the estimates
(\ref{sym-E2.10})-(\ref{sym-E2.13}) and  the Cauchy-Schwarz
inequality, we can obtain (\ref{sym-E2.17}).
\end{proof}

\begin{lem}\label{sym-L2.8}
  Under the assumptions of Lemma \ref{sym-L2.2}, we obtain
        \begin{equation}
          \int^M_0(r^\alpha u_t^2)(x,t)dx+\int^t_0\int^M_0
          \left(r^{2n-2+\alpha}u_{xt}^2+r^{\alpha-2}u_t^2
          \right)dxds\leq C_{19}\epsilon_0^2,\label{sym-E2.36}
        \end{equation}
            \begin{equation}
             \left\|\frac{u}{r}(\cdot,t)
             \right\|_{L^\infty}+ \| (r^{n-1}u)_x(\cdot,t)\|_{L^\infty}\leq
              C_{20}\epsilon_0,\label{sym-E2.37}
            \end{equation}
            \begin{equation}
              \int^M_0\left(
              \frac{u^2}{r^2}+r^{2n-2}u_x^2
              \right)(x,t)dx+\int^t_0\int^M_0u_t^2(x,s)dxds\leq
              C_{21}\epsilon_0^2,\label{sym-E2.17-1}
            \end{equation}
for all $t\in[0,T]$, where $\alpha=\frac{3}{2}-n$.
\end{lem}
\begin{proof}
  We differentiate the equation (\ref{sym-E})$_2$ with respect to
  $t$, multiply it by $u_tr^\alpha$ and integrate it over
  $[0,M]$, using the boundary conditions
  (\ref{sym-Efixbd})-(\ref{sym-Efd}), then derive
    \begin{eqnarray}
      &&\frac{d}{dt}\int^M_0\frac{1}{2}r^\alpha u_t^2dx-\frac{\alpha}{2}\int^M_0r^{\alpha-1}u
      u_t^2dx\nonumber\\
               &=&-\int^M_0\left[(2c_1+c_2)\rho^{1+\theta} (r^{n-1}u)_x
                -A\rho^\gamma+P_\infty
               -2c_1(n-1)\rho^\theta\frac{u}{r}\right]\nonumber\\
               &&\times ((n-1)r^{n-2+\alpha}uu_t)_xdx
    -\int^M_0\partial_t\left[(2c_1+c_2)\rho^{1+\theta} (r^{n-1}u)_x
    -A\rho^\gamma\right.\nonumber\\
    &&\left.+A\rho_\infty^\gamma
    -2c_1(n-1)\rho^\theta\frac{u}{r}\right]
     (r^{n-1+\alpha}u_t)_xdx\nonumber\\
            &&+\int^M_02c_1(n-1)\partial_t(r^{n-1}\rho^\theta\partial_x(\frac{u}{r}))
            r^\alpha u_tdxds-\int^M_0f_tr^\alpha u_tdx\nonumber\\
                &&-\left[(n-1)\Delta  P(r^{n-2+\alpha}uu_t)(M,t)+(\Delta
                P)'(r^{n-1+\alpha}u_t)(M,t) \right]\nonumber\\
    &:=&J_1+J_2+J_3+J_4+J_5.\label{sym-E2.38}
    \end{eqnarray}
  From (\ref{sym-E2.9}), (\ref{sym-E2.25}) and (\ref{sym-E2.27}),
  using the Cauchy-Schwarz inequality, we obtain
            \begin{eqnarray}
              J_1
               &\leq& \epsilon \int^M_0 \left(r^{\alpha-2}u_t^2+
               r^{2n-2+\alpha}u_{xt}^2\right)dx \nonumber\\
        &&+C_\epsilon(1+
        \| (r^{n-1}u)_x\|_{L^\infty_{x}}^2)\int^M_0 \left[
        r^{2n-2+\alpha}u_x^2+r^{\alpha-2}u^2+r^{\alpha}(\rho-\rho_\infty)^2
        \right]dx.
            \end{eqnarray}
From (\ref{sym-E2.9})-(\ref{sym-E2.10}),
  using the same argument in the proof of (\ref{sym-E2.30})
   and the Cauchy-Schwarz inequality, we get
        \begin{eqnarray}
                &&J_2+J_3\nonumber\\
          &=&-\int^M_0\left[(2c_1+c_2)\rho^{1+\theta}(r^{n-1}u_t)_x
          (r^{n-1+\alpha}u_t)_x-2c_1(n-1)\rho^\theta(r^{n-2+\alpha}u_t^2)_x
          \right]dx\nonumber\\
                &&+\int^M_0\big\{(2c_1+c_2)(1+\theta)\rho^{\theta+2}
                [(r^{n-1}u)_x]^2-(n-1)(2c_1+c_2)\rho^{1+\theta}(r^{n-2}u^2)_x
                \nonumber\\
          &&-\gamma\rho^{\gamma+1} (r^{n-1}u)_x-2c_1(n-1)\theta\rho^{\theta+1}
           (r^{n-1}u)_x\frac{u}{r}-2c_1(n-1)\rho^\theta\frac{u^2}{r^2}
                \big\}\nonumber\\
                &&\times\left[(n-1+\alpha)\frac{r^{\alpha-1}u_t}{\rho}
                +r^{n-1+\alpha}u_{tx}\right]dx\nonumber\\
                        &&+2c_1(n-1)\int^M_0\left\{(n-1)r^{n-2+\alpha}u\rho^\theta\left(
                        \frac{u}{r}\right)_xu_t\right.  \nonumber\\
            &&\left.-\theta
  r^{n-1+\alpha}\rho^{\theta+1} (r^{n-1}u)_x \left(\frac{u}{r}\right)_xu_t
  -r^{n-1+\alpha}\rho^\theta \left(\frac{u^2}{r^2}\right)_xu_t\right\}dx
            \nonumber\\
  &\leq&
  -(C_{22}-\epsilon)\int^M_0(r^{2n-2+\alpha}u_{xt}^2+r^{\alpha-2}u_t^2)dx
\nonumber\\
            && +C_\epsilon (1+
            \| (r^{n-1}u)_x\|_{L^\infty_{x}}^2)\int^M_0 \left[
            r^{2n-2+\alpha}u_x^2+r^{\alpha-2}u^2+r^{\alpha}(\rho-\rho_\infty)^2
            \right]dx,
        \end{eqnarray}
    \begin{eqnarray}
      J_4&\leq &\epsilon\int^M_0r^{\alpha-2}u_t^2dx
      +C_\epsilon\int^M_0((1-n)Gxr^{-n}u+\partial_r\Delta fu+\partial_t\Delta f)^2r^{2+\alpha}dx\nonumber\\
                &\leq &\epsilon\int^M_0r^{\alpha-2}u_t^2dxds
              +C_\epsilon(f_2^2+\int^M_0r^{2+\alpha}u^2dx )
    \end{eqnarray}
and
                \begin{eqnarray}
      J_5&=&-\int^M_0\left[(\Delta
      P)'(r^{n-1+\alpha}u_t)_x+(n-1)\Delta P(r^{n-2+\alpha}uu_t)_x
      \right]dx\nonumber\\
                &\leq &\epsilon\int^M_0(r^{\alpha-2}u_t^2+r^{2n-2+\alpha}u_{xt}^2)dx
              +C_\epsilon(|\Delta P|^2+|(\Delta P)'|^2). \label{sym-E2.40}
            \end{eqnarray}
Let $\epsilon=\frac{1}{8}C_{22}$, from
(\ref{sym-E2.38})-(\ref{sym-E2.40}), we have
         \begin{eqnarray*}
          &&\frac{d}{dt}\int^M_0r^\alpha u_t^2  dx+\frac{C_{22}}{2}\int^M_0
          \left(r^{2n-2+\alpha}u_{xt}^2+r^{\alpha-2}u_t^2
          \right)dx\\
                &\leq& C(1+\|(r^{n-1}u)_x\|_{L^\infty_{x}}^2)\int^M_0 \left[
            r^{2n-2+\alpha}u_x^2+r^{\alpha-2}u^2+r^{\alpha}(\rho-\rho_\infty)^2
            \right]dx\nonumber\\
                      &&+C\left(f_2^2
                    +|\Delta P|^2+|(\Delta P)'|^2
                    \right)+C_{23}\epsilon_1 \frac{C_{22}}{2}\int^M_0
                  r^{\alpha-2}u_t^2dx.
        \end{eqnarray*}
If $C_{23}\epsilon_1\leq \frac{1}{2}$, from  (\ref{sym-E2.25}) and
(\ref{sym-E2.27}), we can obtain
    \begin{eqnarray}
          &&\frac{d}{dt}\int^M_0r^\alpha u_t^2dx+\frac{C_{22}}{4}\int^M_0
          \left(r^{2n-2+\alpha}u_{xt}^2+r^{\alpha-2}u_t^2
          \right)dx\nonumber\\
                &\leq& C(1+\|(r^{n-1}u)_x\|_{L^\infty_{x}}^2)\int^M_0 \left[
            r^{2n-2+\alpha}u_x^2+r^{\alpha-2}u^2+r^{\alpha}(\rho-\rho_\infty)^2
            \right]dx\nonumber\\
                      &&+C\left(f_2^2
                    +|\Delta P|^2+|(\Delta P)'|^2
                    \right)\label{sym-E2.76}
        \end{eqnarray}

and
    \begin{eqnarray}
    &&\int^M_0(r^\alpha u_t^2)(x,t)dx+\int^t_0\int^M_0
          \left(r^{2n-2+\alpha}u_{xt}^2+r^{\alpha-2}u_t^2
          \right)(x,s)dxds\nonumber\\
          &\leq&
          C\epsilon_0^2(1+\|(r^{n-1}u)_x\|_{L^\infty_{tx}}^2).
          \label{sym-E2.68}
    \end{eqnarray}

From the equation (\ref{sym-E})$_2$, we have
    $$
    (2c_1+c_2)r^{n-1}\partial_x(\rho^{1+\theta}\partial_x(r^{n-1}u))=u_t+Ar^{n-1}(\rho^\gamma)_x
    +2c_1(n-1)r^{n-2}u(\rho^\theta)_x+f,
    $$
and using the estimates (\ref{sym-E2.7}),
(\ref{sym-E2.9})-(\ref{sym-E2.10}),
(\ref{sym-E2.25})-(\ref{sym-E2.32}) and (\ref{sym-E2.68}),
conclude that
    $$
      \int^M_0r^{2n-2+\alpha}[\partial_x(\rho^{1+\theta}\partial_x(r^{n-1}u))]^2dx\leq
      C\epsilon_0^2(1+\|(r^{n-1}u)_x\|_{L^\infty_{tx}}^2),
    $$
and
    \begin{equation}
        \int^M_0|\partial_x(\rho^{1+\theta}\partial_x(r^{n-1}u))|dx\leq
      C\epsilon_0(1+\|(r^{n-1}u)_x\|_{L^\infty_{tx}}),\label{sym-E2.42}
    \end{equation}
for all $t\in[0,T]$. From (\ref{sym-E2.9}), (\ref{sym-E2.17}) and
(\ref{sym-E2.42}), using Sobolev's embedding Theorem
$W^{1,1}\hookrightarrow L^\infty$, we can obtain
    \begin{equation}
    \|\partial_x(r^{n-1}u)\|_{L^\infty_{tx}}\leq
              C_{24}\epsilon_0(1+\|(r^{n-1}u)_x\|_{L^\infty_{tx}}).
              \label{sym-E2.70}
    \end{equation}
If $C_{24}\epsilon_0\leq \frac{1}{2}$, from (\ref{sym-E2.17}),
(\ref{sym-E2.68}) and (\ref{sym-E2.70}), we can get
(\ref{sym-E2.36})-(\ref{sym-E2.17-1}) immediately.
\end{proof}

Now, we can let
    \begin{equation}
      \epsilon_1=(C_{16}+C_{20}  )\epsilon_0.\label{sym-self-E3.76}
    \end{equation}
If
$(1+\frac{4}{\underline{\rho}}+\frac{C_2}{\delta_5}+\frac{C_4}{C_5}
+4C_8+2C_{14}^2+4C_{17}+2C_{23})\epsilon_1+2C_{24}\epsilon_0 \leq
1$, using the results in Lemmas \ref{sym-L2.2}-\ref{sym-L2.8}, we
finish the proof of the \textbf{Claim 1}. Thus, we can let
$\epsilon_0$ be a positive constant satisfying
    \begin{equation}
(1+\frac{4}{\underline{\rho}}+\frac{C_2}{\delta_5}+\frac{C_4}{C_5}
+4C_8+2C_{14}^2+4C_{17}+2C_{23})(C_{16}+C_{20}
)\epsilon_0+2C_{24}\epsilon_0 = 1.\label{sym-self-E3.77}
        \end{equation}
 Using the classical continuity method, we can obtain the
following lemma.

\begin{lem}\label{sym-L2.9}
Under the assumptions in Theorem \ref{sym-thm}, the solution
$(\rho,u)$ satisfies the estimates
(\ref{sym-E2.9})-(\ref{sym-E2.13}), (\ref{sym-E2.21}),
(\ref{sym-E2.14}), (\ref{sym-E2.24})-(\ref{sym-E2.32-1-1}),
(\ref{sym-E2.36})-(\ref{sym-E2.17-1}) for all $t\geq0$.
\end{lem}

\begin{rem}
  In this section, we just
give some ideas. In deed, we use the continuity method to obtain
the bound of $\frac{u_i(t)}{r_i(t)}$, $i=0,\ldots,N+1$, in Section
\ref{sym-Sec3}. The basic theory of ordinary differential
equations guarantees $\rho_i(t),u_i(t),r_i(t)\in C([0,T^*))$,
$i=0,\ldots,N$. Since $r_j^n(t)\geq h$, $j=0,\ldots,N$, we have
     $\frac{u_i(t)}{r_i(t)}\in C([0,T^*))$, $i=0,\ldots,N$.
Thus, we can use the continuity method.
\end{rem}

From Lemma \ref{sym-L2.9}, we can obtain the following lemma
easily.
\begin{lem}\label{sym-L2.10}
Under the assumptions in Theorem \ref{sym-thm}, if $\epsilon_0$ is
small enough, we have
    $$
      \|\rho(\cdot,t_1)-\rho(\cdot,t_2)\|_{L^2}\leq C|t_1-t_2|,
    $$
        $$
      \|u(\cdot,t_1)-u(\cdot,t_2)\|_{L^2}\leq C|t_1-t_2|,
        $$
    $$
            \|r(\cdot,t_1)-r(\cdot,t_2)\|_{L^\infty}\leq C|t_1-t_2|,
    $$
        $$
                \|\partial_x(r^{n-1}u)(\cdot,t_1)-
                \partial_x(r^{n-1}u)(\cdot,t_2)\|_{L^2}\leq C|t_1-t_2|^{\frac{1}{2}},
        $$
        $$
          \|((r^{n-2}u)_x,(r^{n-1})_x)(\cdot,t)\|_{L^{n-\frac{1}{2}}}\leq C,
        $$
for all $t_1,t_2,t\geq0$.
\end{lem}

\section{Difference scheme and approximate solutions.}\label{sym-Sec3}\setcounter{equation}{0}
In this section, applying a discrete difference scheme as in
\cite{Chen2002}, we construct approximate solutions to the initial
boundary value problem (\ref{sym-E})-(\ref{sym-Efd}).

For any given positive integer $N$,  let $h=\frac{1}{N}$ be an
increment in $x$ and $x_j=jh$ for $j\in\{0,\ldots,N\}$. For each
    integer $N$, we construct the following time-dependent
    functions:
        $$
          (\rho_j(t),u_j(t),r_j(t)),\ j=0,\ldots,N,
        $$
    that form a discrete approximation to $(\rho,u,r)(x_j,t)$
    for $j=0,\ldots,N$.

        First, $\rho_i(t)$, $u_j(t)$ and $r_{i+1}(t)$, $i=0,\ldots,N$, $j=1,\ldots,N$,
        are determined by the following system of $3N+2$
        differential equations:
            \begin{equation}
              \frac{d}{dt}\rho_i=-\rho^2_i\delta(r^{n-1}_iu_i),\label{De1}
            \end{equation}
        \begin{equation}
          \frac{d}{dt}u_j=r^{n-1}_j\delta
          \sigma_j-2(n-1)r^{n-2}_ju_j\delta(\mu_{j-1})-f_j,
          \label{De2}
        \end{equation}
            \begin{equation}
              \frac{d}{dt}r_{i+1}=u_{i+1},
               \label{De3}
            \end{equation}
with the boundary conditions:
    \begin{equation}
      u_0(t)=0,\ r_0^n(t)=h,
      \label{E3.10}
    \end{equation}
        \begin{equation}
          P_N-\rho_N(\lambda_N+2\mu_N)
          \delta(r^{n-1}_Nu_N)+2(n-1)\frac{u_{N+1}}{r_{N+1}}\mu_N=P_\Gamma,
          \label{DEbd}
        \end{equation}
and initial data
    \begin{equation}
      (\rho_j, u_j)(0)=\left(\frac{1}{h}\int^{jh}_{(j-1)h}\rho_0(y)dy,
      \frac{1}{h}\int^{jh}_{(j-1)h}u_0(y)dy\right),j=1,\ldots,N,
    \end{equation}
        \begin{equation}
          \rho_0(0)=\rho_1(0),
        u_0(0)=0,\ r_0^n(0)=h,
        \end{equation}
            \begin{equation}
            r^n_i(0)=h+n\sum^{i-1}_{l=0}\frac{h}{\rho_l(0)},
          \ i=1,\ldots,N+1,
            \end{equation}
and $u_{N+1}(0)$ satisfies
        \begin{equation}
          P_N(0)-\rho_N(0)(\lambda_N(0)+2\mu_N(0))
          \delta(r^{n-1}_N(0)u_N(0))+2(n-1)\frac{u_{N+1}(0)}{r_{N+1}(0)}\mu_N(0)=P_\Gamma(0),
          \label{sym-E3.9}
        \end{equation}
where $\delta$ is the operator defined by $\delta
w_j=(w_{j+1}-w_j)/h$, and
        $$
      \sigma_j(t)=\rho_{j-1}(\lambda_{j-1}+2\mu_{j-1})
      \delta(r^{n-1}_{j-1}u_{j-1})-P_{j-1},
        $$
        $$
          \lambda_{j}=\lambda(\rho_{j}),
          \mu_{j}=\mu(\rho_{j}),
          P_{j}=P(\rho_{j}),
        $$
            $$
              f_j(t)=f(jh,r_j,t).
            $$
The boundary conditions (\ref{E3.10})-(\ref{DEbd}) are consistent
with the initial data. The condition (\ref{DEbd}) determines
$u_{N+1}(t)$.

Let $(\rho_{\infty,i},r^n_{\infty,
i})=(\rho_\infty(ih),h+r^n_\infty(ih))$, $i=0,\ldots,N$, we have
            $$
              r_{\infty, j}^{n-1}\delta (A\rho_{\infty,
              j-1}^\gamma)=-\frac{Gjh}{r_{\infty,j}^{n-1}}+Q_{1j},
            $$
                $$
                  r^n_{\infty, j}=h+n\sum^{j-1}_{k=0}\frac{h}{\rho_{\infty, k}}+Q_{2j},
                $$
and
    $$
    |Q_{1j}|\leq C(jh)^{\frac{1-n}{n}}h,
    \ |Q_{2j}|\leq C(jh)^{\frac{2}{n}}h, \ j=1,\ldots,N.
    $$

         Then, for any small $h$, the initial data $(\rho_i,
u_i,r_i)(0)$ and the external force $f_i$, $i=0,\ldots,N$,
satisfies
    \begin{equation}
      \max_{i\in\{0,\ldots,N\}}|\rho_i(0)-\rho_{\infty,i}|^2+
        \sum^{N-1}_{j=0}\left[
        r_j^{2n-2+\alpha}(\delta \rho_j-\delta \rho_{\infty,j})^2\right](0)h\leq C\epsilon_0^2,
      \label{sym-E3.11}
    \end{equation}
        \begin{equation}
         C^{-1}(i+1)h\leq r^n_i(0)\leq C(i+1)h,\ \sum^{N}_{j=0}\left[
        r_j^{-2}u_j^2+r_j^{2n-2}(\delta u_j)^2\right](0)h\leq C\epsilon_0^2,
        \end{equation}
            \begin{equation}
        \sum^{N}_{j=1}\left(
        r_j^{2n-2\alpha}\left[\delta\left(\rho_{j-1}^{1+\theta}\delta(r_{j-1}^{n-1}u_{j-1})
        \right)        \right]^2
        \right)(0)h\leq C\epsilon_0^2,
        \label{sym-E3.14}
        \end{equation}
where $C>0$ are independent of $h$.

The basic theory of differential equations guarantees the local
existence of smooth solutions $(\rho_i,u_i,r_i)$ ($i=0,\ldots,N$)
to the Cauchy problem (\ref{De1})-(\ref{sym-E3.9}) on an interval
$[0,T^h)$, such that
    $$
    0<\rho_i(t)<\infty, |u_i(t)|,|r_i(t)|<\infty,i=0,\ldots,N,
    $$
with the aid of (\ref{sym-E3.11})-(\ref{sym-E3.14}).

For any fixed $T>0$, by  virtue of Lemma
\ref{sym-L2.1}-\ref{sym-L2.10} and using  similar arguments as in
\cite{Chen2002,hoff92}, we can obtain the following lemma and
prove that the Cauchy problem (\ref{De1})-(\ref{sym-E3.9}) has a
unique solution for $t\in [0,T]$ when $h\leq h_{T,\epsilon_0}$,
where $h_{T,\epsilon_0}>0$ is a constant dependent on $T$ and
$\epsilon_0$.

\begin{lem}\label{sym-L3.1}
For any $h\in(0,h_{T,\epsilon_0}]$,  there exist a positive
constant $C$ independent of $h$ such that
     $$
          \rho_i(t)\in
          \left[\frac{1}{2}\underline{\rho},\frac{3}{2}\bar{\rho}
          \right], \ \|\rho_i(t)-\rho_{\infty,i}\|_{L^\infty}\leq
          C\epsilon_0,
        $$
  $$
    r_l^n(t)\in [C^{-1}(l+1)h,C(l+1)h],
  $$
        $$
       \sum_{j=0}^N\left(u_j^2(t)+ |\rho_j(t)-\rho_{\infty,j}|^2\right)h\leq C\epsilon_0^2,
        $$
  $$
    \left|\frac{u_l(t)}{r_l(t)}\right|\leq C\epsilon_0,
  $$
        $$
        \int^t_0\sum_{j=0}^N\left(u_j^2
        +r_j^{2n-2}(\delta u_{j})^2+\frac{u_j^2}{r_j^2}
        \right)(s)hds\leq C\epsilon_0^2,
        $$
            $$
              \sum_{j=0}^N\left(
              \frac{u_j^2}{r_j^2}+r_j^{2n-2}(\delta u_j)^2
              \right)(t)h+\int^t_0\sum_{j=0}^N\left(\frac{d}{dt}u_j\right)^2(s)hds\leq
              C\epsilon_0^2,
            $$
        $$
          \int^t_0\sum_{j=0}^N\left[r_j^\alpha(\rho_j-\rho_{\infty,j})^2
          +(r_j-r_{\infty,j})^2\right]hds\leq
          C\epsilon_0^2,
        $$
        $$
         \sum_{j=0}^N(r_j^\alpha
         u_j^2)(t)dx+\int^t_0\sum_{j=0}^Nr_j^\alpha\left(
         r_j^{2n-2}(\delta u_j)^2+\frac{u_j^2}{r_j^2}
         \right)hds\leq C\epsilon_0^2,
        $$
        $$
          \sum_{j=0}^{N-1}(r_j^{2n-2+\alpha}(\delta
          \rho_j-\delta\rho_{\infty,j})^2)(t)h+\int^t_0\sum_{j=0}^{N-1}
          (r_j^{2n-2+\alpha}(\delta \rho_j-\delta\rho_{\infty,j})^2)(s)hds\leq C
          \epsilon_0^2,
        $$
                $$
                  \sum_{j=0}^{N-1}|\delta \rho_j-\delta\rho_{\infty,j}|(t)h\leq C
                  \epsilon_0,
                $$
    $$
          \sum_{j=0}^N\left[r_j^\alpha
          \left(\frac{d}{dt}u_j\right)^2\right](t)h\leq C\epsilon_0^2,
        $$
            $$
              |\delta (r_i^{n-1}u_i)(t)|\leq
              C\epsilon_0,
            $$
    $$
     \sum_{j=0}^N \left(|\rho_j(t_1)-\rho_j(t_2)|^2+
     |u_j(t_1)-u_j(t_2)|^2\right)h\leq C|t_1-t_2|^2,
    $$
    $$
            |r_l(t_1)-r_l(t_2)|\leq C|t_1-t_2|,
    $$
        $$
               \sum_{j=0}^N |\delta(r_j^{n-1}u_j)(t_1)-
                \delta (r_j^{n-1}u_j)(t_2)|^2h\leq C|t_1-t_2|,
        $$
        $$
         \sum_{j=0}^N\left( |\delta (r_j^{n-2}u_j)|^{n-\frac{1}{2}}
         +|\delta (r_j^{n-1})|^{n-\frac{1}{2}}\right)h\leq C,
        $$
for all $t_1,t_2,t\in[0,T]$, $i\in\{0,\ldots,N\}$ and
$l\in\{1,\ldots,N+1\}$.
\end{lem}
Now, we can define our approximate solutions
$(\rho^N,u^N,r^N)(x,t)$ for the Cauchy problem
(\ref{sym-E})-(\ref{sym-Efd}). For each fixed $N$ and $t\in[0,T]$,
we define piecewise linear continuous functions
$(\rho^N,u^N,r^N)(x,t)$ with respect to $x$ as follows: when
$x\in[[xN],[xN]+1]$
    $$
\rho^N(x,t)=\rho_{[xN]}(t)+(xN-[xN])(\rho_{[xN]+1}(t)-\rho_{[xN]}(t)),
    $$
    $$
    u^N(x,t)=u_{[xN]}(t)+(xN-[xN])(u_{[xN]+1}(t)-u_{[xN]}(t)),
    $$
    $$
r^N(x,t)=\left(r^n_{[xN]}(t)+(xN-[xN])(r^n_{[xN]+1}(t)-r^n_{[xN]}(t))
\right)^{1/n}.
    $$

From Lemma \ref{sym-L3.1}, using similar arguments as in
\cite{Chen2002,hoff92}, we can obtain the compactness of
approximate solutions $(\rho^N,u^N,r^N)$ and prove the existence
part of Theorem \ref{sym-thm}. Since the constant $C$ in Lemma
\ref{sym-L3.1} is independent of $T$, we can obtain the regularity
estimates (\ref{sym-E1.32})-(\ref{sym-E1.40-1}) easily.

\section{Uniqueness.}\label{sym-Sec4}\setcounter{equation}{0}
In this section, applying energy method, we will prove the
uniqueness of the solution in Theorem  \ref{sym-thm}. Let
$(\rho_1,u_1,r_1)(x,t)$ and $(\rho_2,u_2,r_2)(x,t)$ be two
solutions in Theorem \ref{sym-thm}. Then, we have, $i=1,2$,
$(x,t)\in[0,M]\times[0,T]$
    \begin{equation}
      \rho_i(x,t)\in
      \left[\frac{1}{2}\underline{\rho},\frac{3}{2}\bar{\rho}
          \right],
      \ C^{-1}x^\frac{1}{n}\leq r_i(x,t)\leq Cx^\frac{1}{n},
      \label{sym-E5.2.0}
    \end{equation}
        \begin{equation}
       |x^{-\frac{1}{n}}u_i(x,t)|+|x^{\frac{n-1}{n}}\partial_xu_i(x,t)|\leq C.\label{sym-E5.2.1}
        \end{equation}

For simplicity, we may assume that $(\rho_1,u_1,r_1)(x,t)$ and
$(\rho_2,u_2,r_2)(x,t)$ are suitably smooth since the following
estimates are valid for the solutions with the regularity
indicated in Theorem \ref{sym-thm} by using the Friedrichs
mollifier.

Let
    $$
      \varrho=\rho_1-\rho_2,
      \ w=u_1-u_2,
      \ R=r_1-r_2.
    $$
From (\ref{sym-E2.1-00}), we have
    \begin{eqnarray}
      &&\frac{d}{dt}\int^M_0x^{-\frac{2}{n}}R^2(x,t)dt
            =2\int^M_0x^{-\frac{2}{n}}R wdx\nonumber\\
      &\leq&\epsilon\int^M_0x^{-\frac{2}{n}}w^2dx+C_\epsilon\int^M_0x^{-\frac{2}{n}}R^2dx.\label{sym-E4.4}
    \end{eqnarray}
From (\ref{sym-E}) and (\ref{sym-E5.2.0})-(\ref{sym-E5.2.1}), we
have
    \begin{eqnarray}
    &&\frac{d}{dt}\int^M_0\varrho^2(x,t)dt
                =2\int^M_0\varrho\partial_t(\rho_1-\rho_2)dx\nonumber\\
        &=&2\int^M_0\varrho\left(-\rho_1^2r_1^{n-1}\partial_xu_1+\rho_2^2r_2^{n-1}\partial_xu_2
        -(n-1)\frac{\rho_1u_1}{r_1}+(n-1)\frac{\rho_2u_2}{r_2}\right)dx\nonumber\\
        &\leq &\epsilon\int^M_0(x^{\frac{2n-2}{n}}w_x^2+x^{-\frac{2}{n}}w^2)dx
        +C_\epsilon\int^M_0(\varrho^2+x^{-\frac{2}{n}}R^2)dx.\label{sym-E4.5}
    \end{eqnarray}
From the equation (\ref{sym-E})$_2$ and boundary conditions
(\ref{sym-Efixbd})-(\ref{sym-Efd}), we get
\begin{eqnarray}
      &&\frac{d}{dt}\int^M_0\frac{1}{2}w^2(x,t)dx\nonumber\\
            &&+\int^M_0 \{ (\frac{2}{n}c_1+c_2)\rho_1^{1+\theta} [(r_1^{n-1}w)_x]^2+\frac{2(n-1)}{n}
      c_1\rho_1^{1+\theta}(r_1^{n-1}w_x-\frac{w}{r_1\rho_1})^2 \}dx\nonumber\\
            &=&-\int^M_0\partial_x(r_1^{n-1}w)\left[(2c_1+c_2)(\rho_1^{1+\theta}-
            \rho_2^{1+\theta})\partial_x(r_1^{n-1}u_2)\right.\nonumber\\
      &&\left.+(2c_1+c_2)
            \rho_2^{1+\theta}\partial_x((r_1^{n-1}-r_2^{n-1})u_2)
            -(\rho_1^\gamma-\rho_2^\gamma)\right]dx\nonumber\\
      &&+\int^M_02c_1(n-1)\partial_x\left[
      r_1^{n-1}wu_2(\frac{1}{r_1}-\frac{1}{r_2})\right]\rho_1^\theta
      dx\nonumber\\
            &&+\int^M_02c_1(n-1)\partial_x\left[
            r_1^{n-1}w(\frac{u_2}{r_2})\right](\rho_1^\theta-\rho_2^\theta)dx\nonumber\\
                    &&            +\int^M_02c_1(n-1)\rho_2^\theta\partial_x\left[
                (r_1^{n-1}-r_2^{n-1})w\frac{u_2}{r_2}
                \right]dx\nonumber\\
       &&-\int^M_0\partial_x((r_1^{n-1}-r_2^{n-2})w)\left[(2c_1+c_2)\rho_2^{1+\theta}\partial_x(r_2^{n-1}u_2)
       -\rho_2^\gamma       \right]dx\nonumber\\
                &&+\int^M_0wGx(r_2^{1-n}-r_1^{1-n})dx+\int^M_0w(\Delta f(x,r_2,t)-\Delta f(x,r_1,t))dx.\label{sym-E4.6}
    \end{eqnarray}
 From (\ref{sym-E5.2.0})-(\ref{sym-E5.2.1}) and
(\ref{sym-E4.6}), we have
    \begin{eqnarray}
    &&\frac{d}{dt}\int^M_0\frac{1}{2}w^2(x,t)dx
      +C_{22}\int^M_0\left\{      x^{\frac{2n-2}{n}}w_x^2+x^{-\frac{2}{n}}w^2\right\}dx\nonumber\\
            &\leq&
            C\int^M_0\left(x^{-\frac{2}{n}}R^2+\varrho^2+w^2
            \right)dx.\label{sym-E4.7}
    \end{eqnarray}

From (\ref{sym-E4.4})-(\ref{sym-E4.5}) and (\ref{sym-E4.7}),
letting $\epsilon=\frac{1}{4}C_{22}$,  we obtain
    $$
      \frac{d}{dt}\int^M_0[w^2+\varrho^2+x^{-\frac{2}{n}}R^2]dx
        \leq C\int^M_0\left(x^{-\frac{2}{n}}R^2+\varrho^2+w^2
            \right)dx.
    $$
Using Gronwall's inequality, we have for any $t\in[0,T]$,
    $$
    \int^M_0[w^2+\varrho^2+x^{-\frac{2}{n}}R^2]dx=0.
    $$
This prove the uniqueness of solution in Theorem \ref{sym-thm}.

\section{Asymptotic behavior}\label{sym-Sec5}\setcounter{equation}{0}
In this section, we consider the asymptotic behavior of the
solution to the free boundary problem
(\ref{sym-E})-(\ref{sym-Efd}). We will show that the solution to
the free boundary problem tends to the stationary solution as
$t\rightarrow+\infty$.

The following lemma is proved in \cite{Straskraba}.
\begin{lem}\label{sym-Lem4.1}
  Suppose that $y\in W^{1,1}_{loc}(\mathbb{R}^+)$ satisfies
    $$
      y=y_1'+y_2,
    $$
and
    $$
    |y_2|\leq \sum_{i=1}^n\alpha_i,
    \ |y'|\leq \sum_{i=1}^n\beta_i,
    \ \textrm{ on }\mathbb{R}^+
    $$
where $y_1\in W^{1,1}_{loc}(\mathbb{R}^+)$, and
$\displaystyle{\lim_{s\rightarrow+\infty}}y_1(s)=0$ and
$\alpha_i,\beta_i\in L^{p_i}(\mathbb{R}^+)$ for some
$p_i\in[1,\infty)$, $i=1,\ldots,m$. Then
$\displaystyle{\lim_{s\rightarrow+\infty}}y(s)=0$.
\end{lem}

\begin{prop}\label{sym-Pro5.1}
  Under the assumptions of Theorem \ref{sym-thm}, the total kinetic
  energy
  $$E(t):=\int^M_0\frac{1}{2}u^2(x,t)dx\rightarrow 0
  \textrm{ as }
  \ t\rightarrow+\infty.$$
\end{prop}
\begin{proof}
From (\ref{sym-E2.13}) and Lemma \ref{sym-L2.9}, we have  $
E(t)\in L^1(\mathbb{R}^+)$. Using the Cauchy-Schwarz inequality,
we obtain
    $$
          |E'(t)|\leq E(t)+\int^M_0u_t^2dx.
    $$
Taking into account the estimate (\ref{sym-E2.17-1}) and Lemma
\ref{sym-L2.9}, applying Lemma \ref{sym-Lem4.1}, we finish the
proof.
   \end{proof}

\begin{prop}
 Under the assumptions of Theorem \ref{sym-thm}, we have
  $$\int^M_0(r-r_\infty)^2(x,t)dx\rightarrow 0
  \textrm{ as }
  \ t\rightarrow+\infty.$$
\end{prop}
\begin{proof}
From (\ref{sym-E2.21}) and Lemma \ref{sym-L2.9}, we have  $
\int^M_0(r-r_\infty)^2(x,t)dx\in L^1(\mathbb{R}^+)$. Using the
Cauchy-Schwarz inequality, we obtain
    $$
      \left| \frac{d}{dt}\int^M_0(r-r_\infty)^2 dx\right|=
      \left|2\int^M_0(r-r_\infty)udx\right|\leq 2E(t)+\int^M_0(r-r_\infty)^2 dx.
    $$
Taking into account the estimate  $ E(t)\in L^1(\mathbb{R}^+)$,
applying Lemma \ref{sym-Lem4.1},  we finish the proof.
   \end{proof}

\begin{prop}\label{sym-Pro5.3}
  Under the assumptions of Theorem \ref{sym-thm}, we have
    \begin{equation}
      \int^M_0(\rho-\rho_\infty)^2(x,t)dx\rightarrow0,
      \label{sym-E5.1}
    \end{equation}
 and       \begin{equation}
          \|(\rho-\rho_\infty)(\cdot,t)\|_{L^q}\rightarrow0,
          \ q\in[1,\infty),\label{sym-E5.2}
        \end{equation}
  as $t\rightarrow+\infty$.
\end{prop}
\begin{proof}
  From (\ref{sym-E2.21}) and Lemma \ref{sym-L2.9},
  we have $\int^M_0(\rho-\rho_\infty)^2
  (x,t)dx\in L^1(\mathbb{R}^+)$. From (\ref{sym-E1.32}),
  using the Cauchy-Schwarz inequality, we obtain
    $$
      \left|\frac{d}{dt}\int^M_0(\rho-
      \rho_\infty)^2dx\right|
            \leq \int^M_0(\rho-
      \rho_\infty)^2dx+C\int^M_0(r^{n-1}u)_x^2dx.
    $$
Taking into account the estimate (\ref{sym-E2.13}) and Lemma
\ref{sym-L2.9}, applying Lemma \ref{sym-Lem4.1}, we obtain
(\ref{sym-E5.1}). From (\ref{sym-E1.32}), (\ref{sym-rhoinf}) and
(\ref{sym-E5.1}), we can obtain (\ref{sym-E5.2}) easily.
   \end{proof}

\begin{prop}\label{sym-Pro5.4}
   Under the assumptions of Theorem \ref{sym-thm}, we have
    $$
    \int^M_0x^{\frac{2n-2+\alpha}{n}}((\rho^\theta)_x-(\rho^\theta_\infty)_x)^2(x,t)dx\rightarrow0,
    \textrm { as } t\rightarrow+\infty.
    $$
\end{prop}
\begin{proof}
From (\ref{sym-E1.32}), (\ref{sym-E2.25}), (\ref{sym-E2.32}) and
Lemma \ref{sym-L2.9}, we have
    $$
    \int^M_0x^{\frac{2n-2+\alpha}{n}}((\rho^\theta)_x-(\rho^\theta_\infty)_x)^2(x,t)dx
    \in
    L^1(\mathbb{R}^+).
    $$
From (\ref{sym-E}), (\ref{sym-E1.32})-(\ref{sym-E1.33-1}) and
(\ref{sym-rhoinf}), using the Cauchy-Schwarz inequality, we have
    \begin{eqnarray*}
      &&\left|\frac{d}{dt}\int^M_0x^{\frac{2n-2+\alpha}{n}}((\rho^\theta)_x-(\rho^\theta_\infty)_x)^2 dx
      \right|\\
            &=&2\theta\left|\int^M_0x^{\frac{2n-2+\alpha}{n}}((\rho^\theta)_x-(\rho^\theta_\infty)_x)
            (\rho^{\theta+1}\partial_x(r^{n-1}u))_xdx     \right|\\
      &=&\frac{2\theta}{2c_1+c_2}\left|\int^M_0x^{\frac{2n-2+\alpha}{n}}((\rho^\theta)_x-(\rho^\theta_\infty)_x)
      \left(\frac{u_t}{r^{n-1}}+A(\rho^\gamma)_x\right.\right.\nonumber\\
      &&\left.\left.+2c_1(n-1)\frac{u(\rho^\theta)_x}{r}+\frac{f(x,r,t)}{r^{n-1}}
      \right)dx     \right|\\
      &\leq&C\int^M_0\left[x^{\frac{2n-2+\alpha}{n}}((\rho^\theta)_x-(\rho^\theta_\infty)_x)^2
      +r^{\alpha}u_t^2+r^\alpha(r-r_\infty)^2\right.\nonumber\\
      &&\left.+r^\alpha(\rho-\rho_\infty)^2+r^{\alpha-2}u^2
      \right]dx+f_1^2.
    \end{eqnarray*}
Taking into account the estimates (\ref{sym-E1.32}),
(\ref{sym-E2.25})-(\ref{sym-E2.32}), (\ref{sym-E2.36}) and Lemma
\ref{sym-L2.9}, applying Lemma \ref{sym-Lem4.1}, we end the proof.
   \end{proof}
From Proposition \ref{sym-Pro5.3}-\ref{sym-Pro5.4}, using
Sobolev's embedding Theorem, we can obtain the following corollary
immediately.
\begin{cor}
Under the assumptions of Theorem \ref{sym-thm}, we have
    $$
    \|\rho(\cdot,t)-\rho_\infty(\cdot)\|_{L^\infty}
    +\|r(\cdot,t)-r_\infty(\cdot)\|_{L^\infty}\rightarrow0,
    \textrm { as } t\rightarrow+\infty.
    $$
\end{cor}

\begin{prop}\label{sym-Pro5.5}
Under the assumptions of Theorem \ref{sym-thm}, we have
    $$
    \int^M_0x^\frac{2n-2+\alpha}{n}u_x^2(x,t)dx\rightarrow0,
    \textrm { as } t\rightarrow+\infty.
    $$
\end{prop}
\begin{proof}
  From the estimates (\ref{sym-E1.32}), (\ref{sym-E2.27}) and Lemma \ref{sym-L2.9}, we have
  $$\int^M_0x^\frac{2n-2+\alpha}{n}u_x^2(x,t)dx\in
  L^1(\mathbb{R}^+).$$
  Using  the
Cauchy-Schwarz inequality, we have
  \begin{eqnarray*}
    &&\left|\frac{d}{dt}\int^M_0x^\frac{2n-2+\alpha}{n}u_x^2 dx\right|
        =\left|2\int^M_0x^\frac{2n-2+\alpha}{n}u_xu_{xt}dx\right|\\
            &\leq&\int^M_0x^\frac{2n-2+\alpha}{n}u_x^2 dx
            +\int^M_0x^\frac{2n-2+\alpha}{n}u_{xt}^2 dx.
  \end{eqnarray*}
Taking into account the estimates (\ref{sym-E1.32}),
(\ref{sym-E2.27}), (\ref{sym-E2.36}) and Lemma \ref{sym-L2.9},
applying Lemma \ref{sym-Lem4.1}, we end the proof.
   \end{proof}
From Proposition \ref{sym-Pro5.1} and \ref{sym-Pro5.5}, using
Sobolev's embedding Theorem, we can obtain the following corollary
immediately.
\begin{cor}
Under the assumptions of Theorem \ref{sym-thm}, we have
    $$
    \|u(\cdot,t)\|_{L^\infty}\rightarrow0,
    \textrm { as } t\rightarrow+\infty.
    $$
\end{cor}
Thus, we finish the  proof of Theorem \ref{sym-thm}.

\section{Stabilization rate estimates}\label{sym-Sec6}
Now we are in position to estimate the stabilization rate. We
first state the following proposition which gives the
stabilization rate estimates in $L^2([0,M])$-norm of the solution.
\begin{prop}\label{sym-Pro6.1}
  Under the assumptions of Theorem \ref{sym-thm2}, we have
    \begin{equation}
    \int^M_0\left(u^2+(\rho-\rho_\infty)^2
    +x^{-2}(r^n-r_\infty^n)^2\right)dx\leq Ce^{-a_1t}
    \label{sym-E6.1}
    \end{equation}
and
    \begin{equation}
  |\rho(M,t)-\rho_\infty(M)|+ \left(\int^M_0r^{2n-2}(\rho-\rho_\infty)_x^2dx\right)^\frac{1}{2}
  +\|r(\cdot,t)-r_\infty(x)\|_{L^2}\leq Ce^{-a_1t},\label{sym-E6.2}
    \end{equation}
for all $t\geq0$, where $a_1$ is a positive constant.
\end{prop}
\begin{proof}
Let
    $$V_{1}=\int^M_0\frac{1}{2}u^2dx+S[V]-S[V_\infty],
     $$
     $$W_1=\int^M_0\left\{
              (r^{n-1}u)_x^2+r^{2n-2}u^2_x+\frac{u^2}{r^2}\right\}dx.
     $$
From (\ref{sym-E1.42}), (\ref{sym-E2.16-1})-(\ref{sym-E2.18-1}),
we have
    \begin{equation}
      V_1'+2C_{31}W_1\leq Cf_1V_1^{\frac{1}{2}}+C|\Delta P|^2\leq
      Ce^{-a_0t}V_1^{\frac{1}{2}}+Ce^{-2a_0t},\label{sym-E6.9}
    \end{equation}
    \begin{eqnarray}
    &&C_{32}^{-1}\int^M_0\left(u^2+(\rho-\rho_\infty)^2
    +x^{-2}(r^n-r_\infty^n)^2\right)dx\nonumber\\
        &\leq& V_1 \leq C_{32}\int^M_0\left(u^2+(\rho-\rho_\infty)^2
    +x^{-2}(r^n-r_\infty^n)^2\right)dx,\label{sym-E6.7}
    \end{eqnarray}
and
    \begin{equation}
      C_{33}\|u(\cdot,t)\|_{L^2}\leq W_1.
    \end{equation}
 From (\ref{sym-E2.33-1}), we have
    \begin{eqnarray}
      &&\int^M_0\left[
        (\rho-\rho_\infty)^2+x^{-2}(r^n-r_\infty^n)^2
        \right]dx\nonumber\\
       &\leq& -C_{38}\frac{d}{dt}\int^M_0\frac{u}{r^{n-1}}\left(\frac{r^n}{n}-\frac{r_\infty^n}{n}
      \right)dx+C_{38}W_1+      Ce^{-2a_0t}.
    \end{eqnarray}
From (\ref{sym-E1.32}), we obtain
    \begin{equation}
    \left|C_{38}\int^M_0\frac{u}{r^{n-1}}\left(\frac{r^n}{n}-\frac{r_\infty^n}{n}
      \right)dx\right|\leq C_{39}\int^M_0\left(u^2+|\rho-\rho_\infty|^2
      \right)dx.\label{sym-E6.13}
    \end{equation}
Let
    $$V_2=V_1+\epsilon C_{38}
\int^M_0\frac{u}{r^{n-1}}\left(\frac{r^n}{n}-\frac{r_\infty^n}{n}
      \right)dx,$$
        $$
        W_2=C_{31}W_1+\epsilon\int^M_0\left[
        (\rho^{\gamma}-\rho_\infty^{\gamma})^2+x^{-2}(r^n-r_\infty^n)^2
        \right]dx,
        $$
where
      $\epsilon=\min\{\frac{C_{31}}{C_{38}},\frac{1}{2C_{32}C_{39}}\}$.
 From (\ref{sym-E6.9}) and (\ref{sym-E6.7})-(\ref{sym-E6.13}), we
have
    \begin{equation}
      V_2'+W_2\leq Ce^{-2a_0t},\label{sym-E6.14}
    \end{equation}
    \begin{eqnarray}
    &&C_{39}^{-1}\int^M_0\left(u^2+(\rho-\rho_\infty)^2
    +x^{-2}(r^n-r_\infty^n)^2\right)dx\nonumber\\
        &\leq& V_2
         \leq C_{39}\int^M_0\left(u^2+(\rho-\rho_\infty)^2
    +x^{-2}(r^n-r_\infty^n)^2\right)dx,
    \end{eqnarray}
and
    \begin{equation}
      C_{40}\int^M_0\left(u^2+(\rho-\rho_\infty)^2
    +x^{-2}(r^n-r_\infty^n)^2\right)dx
    \leq W_2. \label{sym-E6.16}
    \end{equation}
 Thus $V_2$ is a Lyapunov functional.  From (\ref{sym-E1.42}),
we obtain the estimate (\ref{sym-E6.1}). From (\ref{sym-E1.32}),
(\ref{sym-E2.36-4}), (\ref{sym-E2.53-1}) and (\ref{sym-E6.1}), we
can get (\ref{sym-E6.2}) easily.
\end{proof}

\begin{prop}
  Under the assumptions of Theorem \ref{sym-thm2}, we obtain
     \begin{equation}
              \int^M_0\left(
              \frac{u^2}{r^2}+r^{2n-2}u_x^2
              \right)(x,t)dx\leq Ce^{-a_3t},\label{sym-E6.23}
            \end{equation}
for all $t\geq0$, where $a_3$ is a positive constant.
\end{prop}
\begin{proof}
Let
    \begin{eqnarray*}
      V_3  &=&\int^M_0\left\{
      \frac{1}{2}(\frac{2}{n}c_1+c_2)\rho^{1+\theta}[(r^{n-1}u)_x]^2+\frac{(n-1)}{n}
      c_1\rho^\theta(r^{n-1}u_x-\frac{u}{r\rho})^2\right.\nonumber\\
                &&\left.+(A\rho_\infty^\gamma-A\rho^\gamma+\Delta P)(r^{n-1}u)_x+
                u\left(G\frac{x}{r^{n-1}}-G\frac{xr^{n-1}}{r_\infty^{2n-2}}
                \right)\right\}dx.
    \end{eqnarray*}
From (\ref{sym-E1.33-1}) and (\ref{sym-E2.45}), we have
    $$
            V_3'                +\int^M_0\frac{1}{2}u_t^2(x,s)dx
            \leq C_{41}\left(f_1^2+|(\Delta
                 P)'|^2+W_2\right).
    $$
From (\ref{sym-E2.10})-(\ref{sym-E2.11}), we have
    $$
        V_3\geq\int^M_0\{C_{42}((r^{n-1}u)_x^2+\frac{u^2}{r^2}+r^{2n-2}u_x^2
      )-C_{43}((\rho-\rho_\infty)^2+|\Delta
      P|^2)\}dx,
      $$
      $$
         V_3
      \leq
      \int^M_0\{C_{44}((r^{n-1}u)_x^2+\frac{u^2}{r^2}+r^{2n-2}u_x^2
      )+C_{43}((\rho-\rho_\infty)^2+|\Delta
      P|^2)\}dx.
        $$
Letting $V_4=V_2+\eta V_3+\eta C_{43}|\Delta P|^2$, where
$\eta=\min\{\frac{1}{2},
\frac{1}{4C_{39}C_{43}},\frac{1}{2C_{41}}\}$. From
(\ref{sym-E6.14})-(\ref{sym-E6.16}), we have
    $$
  CW_2\geq  V_4\geq C^{-1}W_1,
    $$
and
    $$
    V_4'+C^{-1}W_2\leq C\left(f_1^2+|\Delta P|^2+|(\Delta P)'|^2\right).
    $$
Thus $V_4$ is a Lyapunov functional. From (\ref{sym-E1.42}), we
can obtain the estimate (\ref{sym-E6.23}).
\end{proof}

\begin{prop}\label{sym-self-P7.3}
  Under the assumptions of Theorem \ref{sym-thm2},  we obtain
        \begin{equation}
          \int^M_0r^{\frac{1}{2}-m}(\rho-\rho_\infty)^2dxds\leq
          Ce^{-at},\label{sym-E6.34}
        \end{equation}
         \begin{equation}
          \int^M_0r^{\frac{1}{2}-m}(r-r_\infty)^2dxds\leq
          Ce^{-at},\label{sym-E6.35}
        \end{equation}
        \begin{equation}
         \int^M_0(r^{\frac{1}{2}-m}
         u^2)(x,t)dx\leq Ce^{-at},\label{sym-E6.36}
        \end{equation}
        \begin{equation}
          \int^M_0(r^{2n-2+{\frac{1}{2}-m}}(\rho-\rho_\infty)_x^2)(x,t)dx\leq Ce^{-at},
          \label{sym-E6.37}
        \end{equation}
for all $t\geq0$ and $m=0,1,\ldots,n-1$, where $a$ is a positive
constant.
\end{prop}
\begin{proof}
From (\ref{sym-E6.1})-(\ref{sym-E6.2}), we know that the estimates
(\ref{sym-E6.34})-(\ref{sym-E6.37}) hold with $m=0$.

\textbf{Claim 3}: If that (\ref{sym-E6.34})-(\ref{sym-E6.37}) hold
with $m\leq k$, $k\in [0,n-2]$, then the estimates
(\ref{sym-E6.34})-(\ref{sym-E6.37}) hold with $m= k+1$.

We could prove Claim 3 as follows. Let $\alpha_k=\frac{1}{2}-k-1$.
From (\ref{sym-E2.53}), (\ref{sym-E6.1}) and
(\ref{sym-E6.37})($m=k$), we have
    $$
      \int^M_0r^{\alpha_k}(\rho-\rho_\infty)^2dx
                \leq Ce^{-at},
   $$
and (\ref{sym-E6.34}) ($m=k+1$) holds. From (\ref{sym-E1.32}) and
(\ref{sym-E6.34}), we can obtain (\ref{sym-E6.35}) ($m=k+1$)
easily.

 From (\ref{sym-E2.63}), we obtain
        \begin{equation}
         \frac{d}{dt}\int^M_0r^{\alpha_k}
         u^2  dx+C_{42}\int^M_0r^{\alpha_k}\left(
         r^{2n-2}u_x^2+\frac{u^2}{r^2}
         \right)dx
         \leq C\left(f_1^2+|\Delta P|^2
         \right)+Ce^{-at}.\label{sym-E6.44-1}
        \end{equation}
Thus $\int^M_0(r^{\alpha_k}
         u^2)(x,t)dx$ is a Lyapunov functional, and we obtain
(\ref{sym-E6.36}) ($m=k+1$) immediately.

From (\ref{sym-E2.66}) and (\ref{sym-E6.34})-(\ref{sym-E6.36}), we
have
    \begin{eqnarray*}
    &&\frac{d}{dt}\int^M_0H_1^2 dx+\frac{C_{18}}{2}\int^M_0H_1^2 dx\\
    &\leq&
    C\int^M_0\left(r^{\alpha_k}u^2+r^{\alpha_k}(\rho-\rho_\infty)^2+
    r^{\alpha_k}(r-r_\infty)^2
    \right)+Cf_1^2\leq Ce^{-at}.
    \end{eqnarray*}
Thus $\int^M_0H_1^2(x,t)dx$ is a Lyapunov functional. Using the
estimates (\ref{sym-E1.32}) and
(\ref{sym-E6.34})-(\ref{sym-E6.36}), we obtain (\ref{sym-E6.37})
($m=k+1$), finish the proof of \textbf{Claim 3} and Proposition
\ref{sym-self-P7.3} immediately.
\end{proof}

From (\ref{sym-E1.32}), (\ref{sym-E6.34}) and (\ref{sym-E6.37}),
using H\"{o}lder's inequality and Sobolev's embedding Theorem, we
could obtain the following proposition.
\begin{prop}
  Under the assumptions of Theorem \ref{sym-thm2},  we obtain
        $$
          \|\rho(\cdot,t)-\rho_\infty(\cdot)\|_{L^\infty}+
                    \|r(\cdot,t)-r_\infty(\cdot)\|_{L^\infty}\leq
                    Ce^{-at},
        $$
for all $t\geq0$, where $a$ is a positive constant.
\end{prop}

\begin{prop}
  Under the assumptions of Theorem \ref{sym-thm2},  we obtain
        \begin{equation}
          \int^M_0(r^\alpha u_t^2)(x,t)dx\leq Ce^{-at},
          \label{sym-E6.48}
        \end{equation}
for all $t\geq0$, where $\alpha=\frac{3}{2}-n$ and $a$ is a
positive constant.
\end{prop}
\begin{proof}
From (\ref{sym-E1.33-1}) and (\ref{sym-E2.76}), we have
         \begin{eqnarray}
          &&\frac{d}{dt}\int^M_0 r^\alpha u_t^2 dx+C_{45}\int^M_0
          \left(r^{2n-2+\alpha}u_{xt}^2+r^{\alpha-2}u_t^2
          \right)dx\nonumber\\
                &\leq&C_{46}\int^M_0 [
        r^{2n-2+\alpha}u_x^2+r^{\alpha-2}u^2+r^{\alpha}(\rho-\rho_\infty)^2
        ]dx  \nonumber\\
        &&   +           C\left(f_2^2+|\Delta P|^2+|(\Delta
              P)'|^2\right).\label{sym-E6.54}
        \end{eqnarray}
Let $V_4=\int^M_0(r^{\alpha}
         u^2)(x,t)dx+\frac{C_{42}}{2C_{46}}\int^M_0(r^\alpha
u_t^2)(x,t)dx$. From (\ref{sym-E1.42}), (\ref{sym-E6.34}),
(\ref{sym-E6.44-1})($k=n-2$) and (\ref{sym-E6.54}), we have
    $$
      V_4'+C^{-1}V_4\leq Ce^{-a_0t}+Ce^{-at}.
    $$
Thus $V_4$ is a Lyapunov functional, and  we obtain
(\ref{sym-E6.48}) immediately.
\end{proof}

\begin{prop}\label{sym-Pro6.8}
    Under the assumptions of Theorem \ref{sym-thm2},  we obtain
        $$
          \int^M_0r^{2n-2+\alpha}[\partial_x(\rho^{1+\theta}\partial_x(r^{n-1}u))]^2dx
          + \left\|\left(\frac{u}{r},(r^{n-1}u)_x\right)(\cdot,t)
    \right\|_{L^\infty} \leq
              Ce^{-at},
        $$
for all $t\geq0$, where $\alpha=\frac{3}{2}-n$ and $a$ is a
positive constant.
\end{prop}
\begin{proof}
From the equation (\ref{sym-E})$_2$, we have
    $$
    (2c_1+c_2)r^{n-1}\partial_x(\rho^{1+\theta}\partial_x(r^{n-1}u))=u_t+Ar^{n-1}(\rho^\gamma)_x
    +2c_1(n-1)r^{n-2}u(\rho^\theta)_x+f,
    $$
and using the estimates (\ref{sym-E1.32})-(\ref{sym-E1.33-1}),
(\ref{sym-E1.42}), (\ref{sym-E6.34})-(\ref{sym-E6.37}) and
(\ref{sym-E6.48}), conclude that
    $$
      \int^M_0r^{2n-2+\alpha}[\partial_x(\rho^{1+\theta}\partial_x(r^{n-1}u))]^2dx\leq
      Ce^{-at},
    $$
and
    \begin{equation}
        \int^M_0|\partial_x(\rho^{1+\theta}\partial_x(r^{n-1}u))|dx\leq
      Ce^{-at}.\label{sym-E6.56}
    \end{equation}
From (\ref{sym-E1.32}), (\ref{sym-E6.23}) and (\ref{sym-E6.56}),
using Sobolev's embedding Theorem $W^{1,1}\hookrightarrow
L^\infty$, we can obtain
    $$
    \|\partial_x(r^{n-1}u)(\cdot,t)\|_{L^\infty}\leq
              Ce^{-at},
    $$
and
    $$
    \left\|\frac{u}{r}(\cdot,t)
    \right\|_{L^\infty}\leq Ce^{-at}.
    $$
\end{proof}

\end{document}